\documentclass [reqno,10pt]{article}

\usepackage{caption}
\usepackage{pb-diagram}
\usepackage{mathrsfs}
\newcommand{\cat}{CAT$(\Kl)$}
\newcommand{\alex}{\angle}
\usepackage[nointlimits]{amsmath}
\usepackage{amssymb}
\usepackage{amscd}
\usepackage{amsthm}
\usepackage{graphics} 
\usepackage{floatflt}
\usepackage{wrapfig}
\usepackage{nicefrac}

\newtheorem{teo}{Theorem}
\newtheorem{prop}{Proposition}
\newtheorem{lemma}{Lemma}
\newtheorem{remark}{Remark} 
\newtheorem{cor}{Corollary}
\newtheorem{definiz}{Definition}

\newcommand{\rqua}{r_4}
\newcommand{\rcin}{r_5}
\newcommand{\rsei}{r_6}

\newcommand{\czero}{c_0}
\newcommand{\cuno}{c_1}
\newcommand{\cdue}{c_2}
\newcommand{\ctre}{c_3}

\newcommand{\kl}{\kappa}
\newcommand{\Kl}{\kappa}
\newcommand{\meno}{^{-1}}

\newcommand{\inte}{\operatorname {int}}
\newcommand{\um}{{\nicefrac{1}{m}}}

\newcommand{\Bd}{\mathfrak{B}}
\newcommand{\intT}{ T}
\newcommand{\barT}{\overline{T}}
\newcommand{\sect}{\mathfrak{S}}
\newcommand{\Sec}{\S }
\newcommand{\deltone}{\frac{\pi}{8}}
\newcommand{\deltino}{\frac{\pi}{2}}
\newcommand{\opt}{t^+}
\newcommand{\opmt}{t^-}
\newcommand{\oppmt}{t^\pm}
\newcommand{\T}{\Delta}
\newcommand{\ango}{\sphericalangle}
\newcommand{\angor}{\sphericalangle_{\operatorname{or}}}
\newcommand{\no}{^0}
\newcommand{\cl}{{C}}
\newcommand{\intE}{\operatorname{int}E}
\newcommand{\alfa}{\alpha}

\newcommand{\dalf}{\dot{\alpha}}
\newcommand{\undalf}{\underline{\dalf}}
\newcommand{\unbeta}{\underline{\beta}}
\newcommand{\db}{\dot{\beta}}
\newcommand{\undb}{\underline{\db}}
\newcommand{\ga}{\gamma}
\newcommand{\dga}{\dot{\gamma}}
\newcommand{\unga}{\underline{\gamma}}
\newcommand{\unl}{\underline{L}}
\newcommand{\dunga}{\underline{\dga}}

\newcommand{\ddga}{\ddot{\gamma}}

\newcommand{\cobra}{\mathfrak{F}}
\newcommand{\uncobra}{\underline{\cobra}}
\newcommand{\camel}{\mathfrak{G}}
\newcommand{\marmo}{\operatorname{\mathfrak{b}}}
\newcommand{\solle}{\mathfrak{a}} 
\newcommand{\inj}[1]{\operatorname{inj}_{#1}}
\newcommand{\A}{\mathcal{P}}
\newcommand{\eps}               { \varepsilon}
\renewcommand{\phi}             {\varphi}

\newcommand{\restr}[1]          {\raisebox{-.2ex}{\big |}\phantom{}_{  #1}}
\newcommand{\R}                 { \mathbb { R } }
\newcommand{\C}                 { \mathbb { C } }

\newcommand{\Zeta}              { \mathbb {Z}  }

\newcommand{\vol}               {\operatorname{\mathrm {vol}}}
       
\newcommand{\Keler}             {K\"{a}hler }

\renewcommand{\Re}              {\operatorname {Re}}
\renewcommand{\Im}              {\operatorname{Im}}

\newcommand{\Lo}{L}
\newcommand{\enf}               {\emph}
\newcommand{\om}{\omega}
\newcommand{\ra}                {\rightarrow}

\newcommand{\demi}              { \frac {1} {2}  } 
 
\newcommand{\menuno}            {{-1}}  
\newcommand{\men}            {{-1}}  
\newcommand{\vacuo}             {\emptyset}

\newcommand{\la}{\lambda}
\newcommand{\compsubset}{\subset\!\!\! \subset}
\newcommand{\sef}{{B}}

\newcommand{\z}{z}
\newcommand{\sing}{{\operatorname{sing}}}
\newcommand{\reg}{{\operatorname{reg}}}
\newcommand{\Xr}{{X_{\operatorname{reg}}}}

\newcommand{\cc}{\check{c}}
\newcommand{\ux}{{U_xX}}
\newcommand{\mult}{\operatorname{mult}}

\newcommand{\dom}{{d}}
\newcommand{\Dom}{\mathscr{D}}
\newcommand{\diam}{\operatorname{diam}}

\newcommand{\thext}{\theta_{\operatorname{ext}}}
\newcommand{\thint}{\theta_{\operatorname{int}}}

\newcommand{\bx}{B_x}

\begin{document}
\title{Alexandrov curvature of \Keler curves}
 \author{Alessandro Ghigi \footnote{Partially supported by MIUR PRIN 2005
     {\em ``Spazi di moduli e Teorie di Lie''}.}}
%
%

\maketitle

{\abstract{
    \noindent We study the intrinsic geometry of a one-dimensional
    complex space provided with a \Keler metric in the sense of
    Grauert. We show that if $\Kl$ is an upper bound for the Gaussian
    curvature on the regular locus, then the intrinsic metric has
    curvature $\leq \Kl$ in the sense of Alexandrov.

  }
  
}

\tableofcontents{}

\section{Introduction}

%
%
%

Let $\Omega $ be a domain in $\C^n$ and let $X \subset \Omega$ be an
analytic subset, that is a set of the form $X=\{z\in \Omega:
f_1(z)=\cdots =f_N(z)=0\}$ for some functions $f_j$ holomorphic on
$\Omega$.  Denote by $\langle \ , \, \rangle$ the flat metric on
$\C^n$ and by $g$ the Riemannian metric induced on the regular part
$\Xr$ of $X$.  Define a distance on $X$ by setting $d(x,y)$ equal to
the infimum of the lengths of curves lying in $X$ and joining $x$ to
$y$.  Then $(X,d)$ is an intrinsic metric space.  If $X$ is smooth,
$(X,g)$ is a \Keler manifold and one can study the metric properties
of $(X,d)$ using the methods of Riemannian geometry.  If $X$ contains
singularities it is natural to study $(X,d)$ using the notions and
methods of nonregular Riemannian geometry developed by the great
soviet mathematician A.D. Alexandrov and his school (see e.g.
\cite{aleksandrov-zalgaller}, \cite{rinow},
\cite{aleksandrov-berestovskij-nikolaev-eng},
\cite{berestovskij-nikolaev-enci}, \cite{ballmann-nonp},
\cite{nikolaev-ln-synthetic}, \cite{nikolaev-ln-metric},
\cite{bridson-haefliger-libro}, \cite{burago-burago-ivanov}).

The purpose of this paper is to investigate the Alexandrov geometry of
$(X,d)$ in the simplest case, namely when $\dim_\C X=1$.

More generally, we consider the following situation. 
Let $X$ be a one-dimensional connected reduced complex space and let
$\om$ be a \Keler form on $X$ in the sense of Grauert
\cite{grauert-exzeptionelle}.  This means that $\om$ is a \Keler
metric on the regular part of $X$ and that for any singular point
there is a representation of $X$ as an analytic set in some open set
$\Omega \subset \C^n$ such that $\om$ extends to a smooth \Keler
metric on $\Omega$
 (see \S \ref{distance} for precise definitions).  
The \Keler form $\om$ and the complex structure determine a Riemannian
metric $g$ on $\Xr$ which allows to compute the length of paths in
$X$.  Given two point $x,y \in X$ let $d(x,y)$ be the infimum of the
lengths of paths in $X$ from $x$ to $y$.
It turns out that $d$ is an intrinsic distance on $X$ inducing the
original topology.  We refer to $d$ as the intrinsic distance of
$(X,\om)$.  Our main result is the following.  \setcounter{teo}{7}
\begin{teo}
  \label{main}
  Let $X$ be a one-dimensional connected reduced complex space. Let
  $\om$ be a \Keler metric on $X$ in the sense of Grauert and let $d$
  be the intrinsic distance of $(X,\om)$.  If $\Kl$ is an upper bound
  for the Gaussian curvature of $ g$ on $\Xr$, then $(X,d)$ is a
  metric space of curvature $\leq \Kl$ in the sense of Alexandrov.
\end{teo}

\setcounter{teo}{0}

This result is strictly related to a theorem of Mese, of which the
author was not aware until the completion of this work. We explain
briefly the relation among the two results.  An immediate corollary of
our main result is the following.
\begin{cor} \label{corollo} Let $(M,g)$ be a \Keler manifold with
  sectional curvature $\leq \Kl$ for some $\Kl\in \R$.  Let $X\subset
  M$ be a one-dimensional analytic subset and let $d$ be the intrinsic
  distance on $X$.  Then $(X,d)$ is an inner metric space of curvature
  $\leq \Kl$ in the sense of Alexandrov.
\end{cor}
If $X$ is smooth this follows from Gauss equation together with the
\Keler property of $g$.  If $X$ contains singularities it is enough to
apply the Gauss equation on $\Xr$ and Theorem \ref{main}.  More
generally, if $(M,g)$ is a Riemannian manifold with sectional
curvature $\leq \Kl$ and $X$ is a smooth minimal surface in $M$, Gauss
equation implies that $X$ has curvature $\leq \Kl$ as well.
%
%
This leads to the following problem: if $(Y,d)$ is a metric space with
Alexandrov curvature $\leq \Kl$ and $X\subset Y$ is a minimal surface
(in a sense to be defined) is it true that $X$ with the induced metric
has curvature $\leq \Kl$?  Mese \cite{mese-curvature-singular} has
proved this for surfaces that are conformal and energy-minimizing and
Corollary \ref{corollo} immediately follows from her result.  We also
mention that Petrunin \cite{petrunin-metric-minimizing} has gotten the
same result for \emph{metric minimizing} surfaces.


We wish to stress that in Theorem \ref{main} there are no assumptions
on the curvature of the ambient manifold. In fact the ambient manifold
is not specified at all. Therefore our result is
stronger. Nevertheless it follows immediately from Mese theorem that a
one-dimensional analytic subset with intrinsic metric is CAT$(\Kl')$
for some $\Kl'$ possibly larger than $\Kl$. This immediately yields
uniqueness of geodesics in the small, which is one of the hardest part
in our proof. Nevertheless Mese arguments depends on quite deep
analytic tools, while ours uses only the normalisation of a
one-dimensional singularity, Jordan curve theorem, Rauch theorem and
Klingenberg lemma from Riemannian geometry, and some basic constructs
of Alexandrov geometry. Moreover our analysis yields a very concrete
description of all the geometric objects involved.

Therefore we believe that it is still of some interest to present the
proof of Theorem \ref{main} in this form.

\medskip{}

The plan of the paper is the following.


In \S \ref{distance} we recall the definition of \Keler forms on a
singular space, define the intrinsic distance in the one-dimensional
case and prove some basic properties.  Many statements hold in more
general situations, but we restrict from the beginning to the
one-dimensional case in order not to burden the presentation.  At the
end we show that to investigate local problems one might restrict
consideration to the case in which $X$ is a one-dimensional analytic
subset in $\C^n$ provided with a general \Keler metric.
Appropriate conventions and notations are fixed to be used in the
study of this particular case under the additional hypothesis that
there is only one singularity which is (analytically) irreducible.
This study occupies \Sec\Sec \ref{sec-regu}--\ref{convex-section}.

In \Sec \ref{sec-regu} we consider differentiability properties of
segments $\alfa : [0,L] \ra X$. Since $X\subset \C^n$ we can consider
the tangent vector $\alfa(t)$ at least when $\alfa(t) \in \Xr$. The
main point is a H\"older estimate for $\dalf$ (Theorem
\ref{W-Hoelder}). This is proved by expressing the second fundamental
form of $X\subset \C^n$ in terms of the normalisation map. Here is
where the \Keler property is used. Next we make various observations
regarding the asymptotic behaviour of the distance $d$ and of the
tangents to segments close to a singular point.

In \Sec \ref{regu-delta} we study regularity properties of segments
through the normalisation. This is useful to compute angles between
the tangent vectors at a singular point.

\Sec \ref{sec-uni} is the most technical section.  We study uniqueness
properties of geodesics near the singular point.  We construct a
decreasing sequence of radii $r_1 > r_2>r_3 >\rqua >r_5> \rsei$ such that the
geodesic balls centred at the singular point have better and better
uniqueness properties.  As a first step (Prop. \ref{unicita-regolare}
and Theorem \ref{unicita}) we show that if two segments have the same
endpoints then the singular point lies in the interior of the closed
curve formed by the segments.  To prove this we combine extrinsic and
intrinsic information. The former amounts to the H\"older estimate
alluded to above and the finiteness of the area of $X$ (Lelong
theorem).  The latter is provided by Gauss--Bonnet and Rauch theorems.
The Jordan separation theorem is used on several occasions.  Th next
step (Theorem \ref{unicita-giro}) is to show that if two points
sufficiently close to the singularity are joined by two distinct
segments one of them has to pass through the singular point. Here the
argument is based on
the winding number and the fact that $X$ is a ramified covering of the
disc.

In \Sec \ref{convex-section} we prove that sufficiently small balls
centred at the singular point are geodesically convex
(Cor. \ref{Whitehead}). On the way we prove (using an idea from
\cite{mccann-polar}) that the distance from a singular point is
$\cl^1$ in a (deleted) neighbourhood of it.  We establish various
technical properties of segments emanating from the singular point and
the angle their tangent vectors form at the singular point. In
particular we study ''sectors'' with vertex at the singular point
(Lemma \ref{dentro-dentro}) and establish their convexity (Theorem
\ref{settore-convesso}).

In \Sec \ref{Alex-sec} we recall the main concepts of the intrinsic
geometry of metric spaces in the sense of A.D. Alexandrov.  Next, by
combining the information on sectors and angles collected before, we
show that a sufficiently small ball centred at a singular point is a
\cat--space. This completes the proof of Theorem \ref{main} in the
case of an irreducible singularity.  The case of reducible
singularities is dealt with by reasoning as in Reshetnyak gluing
theorem and invoking the result in the irreducible case.

At the end of the paper we observe that the statement corresponding to
Theorem \ref{main} with lower bounds on curvature instead of upper
bounds is false. In particular a \Keler curve $(X,d)$ can have
curvature bounded below in the sense of Alexandrov only if $X$ is
smooth  (Theorem \ref{no-low}).

\paragraph{Acknowledgements}
The author wishes to thank Prof. Giuseppe Savar\'e for turning his
attention to the Alexandrov notions of curvature and both him and
Prof.  Gian Pietro Pirola for various interesting discussions on
subjects connected with this work.  He is also grateful to an
anonymous referee for pointing out the result of Mese
\cite{mese-curvature-singular}.  He also acknowledges generous support
from MIUR PRIN 2005 ``Spazi di moduli e Teorie di Lie''.

\section{Intrinsic distance}
\label{distance}

A Let $X$ be \enf{complex curve}, that is a one-dimensional reduced
complex space.  By definition for any point $x\in X$ there is an open
neighbourhood $U$ of $x$ in $X$, a domain $\Omega$ in some affine space
$\C^n$ and a map $\tau : U \ra \Omega$ that maps $U$ biholomorphically
onto some one-dimensional analytic subset $A\subset \Omega$. We call
the quadruple $(U,\tau, A, \Omega)$ a \enf{chart} around $x$.
\begin{definiz}
  A \enf{\Keler form} on $X$ is a \Keler form $\om$ on $X_\reg$ with
  the following property: for any $x\in X_\sing$ there is are a chart
  $(U,\tau, A, \Omega)$ around $x$ and a \Keler form $\om'$ on
  $\Omega$ such that $\tau^*\om' = \om$ on $U\cap X_\reg$.  We call
  $\om'$ a \enf{local extension} of $\om$.
\end{definiz}
This definition is due to Grauert \cite[\S
3.3]{grauert-exzeptionelle}.  A \enf{\Keler curve} is a complex curve
with a fixed \Keler form.  Let $(X,\om)$ be a \Keler curve. Denote by
$J$ the complex structure on $\Xr$. Then $g(v,w) = \om(v, Jw)$ defines
a Riemannian metric on $\Xr$.  Denote by $|v|_g$ the norm of $v\in
T_x\Xr$ with respect to $g$.  A path $\alfa : [a,b] \ra X$ is of class
$\cl^1$ if $\tau \circ \alfa$ is $\cl^1$ for any chart.
For a piecewise $\cl^1$ path $\alfa $ the length is defined by
\begin{equation}
  L(\alfa)= \int_{\alfa\meno (\Xr)} |\dalf(t)|_g dt.
\end{equation}
\begin{lemma}\label{vero-length}
  Let $(U,\tau,A, \Omega)$ be a chart and $\om'$ a \Keler form on
  $\Omega$ extending $\om$, with $g'$ the corresponding metric. If
  $\alfa : [a,b]\ra U$ is a piecewise $\cl^1$ path and $\beta=\tau
  \circ \ga$, then
  \begin{equation}
    \label{eq:length=length}
    L(\alfa) = \int_a^b |\dot{\beta}(t)|_{g'} dt
  \end{equation}
\end{lemma}
\begin{proof} Let $E=\alfa^\menuno{}(X_\reg)$, $F=I\setminus E$,
  $B=F^0$, $D=\partial F$. Then $I=E \sqcup B \sqcup D$. Since $X$ has
  isolated singularities $\alfa$ and $\beta$ are constant on the
  connected components of $F$, so $\dot{\beta} \equiv 0$ on $B$. The
  set $D$ is countable, so has zero measure. Therefore
  \begin{equation*}
    \int_a^b |\dot{\beta}|_{g'} dt = 
    \int_E |\dot{\beta}|_{g'} = \int_E |\dalf|_g
    =  L(\alfa) .
  \end{equation*}
\end{proof}
For $x,y \in X$ set
\begin{gather}
  \label{eq:d-om}
\begin{aligned}
  d(x,y)=\inf\{L(\alfa) : \ & \alfa \text{ piecewise } C^1 \text{ path in
  } X \text{ with } \\ &\alfa(0)=x, \, \alfa(1)=y\}.
\end{aligned}
\end{gather}
For $r>0$ set also $ \Bd(x,r) = \{y\in X: d(x,y)<r\}$.  Recall the
following fundamental result of \L ojasiewicz.
\begin{teo}
  [\L ojasiewicz, {\cite [\S 18, Prop. 3, p.97]{lojasiewicz-ihes}}]
  \label{loja} Let $A$ be an analytic subset in a domain $\Omega
  \subset \C^n$ and $z_0\in A$. Then there are $C>0$, $ \mu \in (0,
  1]$ and a neighbourhood $V$ of $z_0$ in $A$ such that for any
  $z,z'\in V$ there is a real analytic path $\beta: [0,1]\ra A$
  joining $z$ to $z'$ with $ \int _0^1 |\dot{\beta}| dt \leq C
  |z-z'|^\mu$. (Here $| \cdot |$ denotes the Euclidean norm in
  $\C^n$.)
\end{teo}
\begin{prop}
  If $(X, \om)$ is a connected \Keler curve, then $d$ is a distance on
  $X$ inducing the original topology.
\end{prop}
\begin{proof}
  We start by showing that $d(x,y) $ is finite for any $(x,y) \in
  X\times X$.  If $x$ and $y$ belong to the same connected component
  of $\Xr$ this is obvious. Assume that $x \in X_{\sing}$. Let
  $(U,\tau,A,\Omega)$ be a chart around $x$ and $\om'$ a local
  extension of $\om$. By restricting $U$ we may assume that there is a
  constant $ C>0$ such that $C\meno |d\tau( v) | \leq |v|_g \leq C
  |d\tau(v)| $ for any $v \in TU_\reg$.  If $\alfa : [a, b] \ra U$ is
  a piecewise $\cl^1$ curve and $\beta = \tau\circ \alfa$, then
  $C\meno L(\beta) \leq L(\alfa) \leq C L(\beta)$, where the length of
  $\beta $ is computed with respect to the Euclidean norm.  By \L
  ojasiewicz Theorem 
  for any point $y\in U$ there is a piecewise $\cl^1$ path $\alfa$ in
  $A$ joining $\tau(x) $ to $\tau(y)$. Then $\beta=\tau\meno \circ
  \alfa$ is a path in $X$ joining $x$ to $y$ with $L(\beta) \leq C
  \cdot L(\alfa) < +\infty$ hence $d(x,y) < + \infty$ for all $y\in
  U$.  Because the length functional $L$ is additive with respect to
  the concatenation of paths, it follows that $d(x,y) < + \infty $ for
  all $y$ in some irreducible component of $X$ that passes through
  $x$. Since $X$ is connected this yields finiteness of $d$.

  At this point one might apply the general machinery of
  \cite[p.123ff]{rinow} or \cite[p.26ff]{burago-burago-ivanov}.
%
%
  The class of piecewise $\cl^1$ paths is closed under restriction,
  concatenation and $\cl^1$ reparametrisations.  Moreover $L$ is
  invariant under $\cl^1$ reparametrisation, it is an additive
  function on the intervals and $L(\alfa\restr{[a,t]})$ is a
  continuous function of $t \in [a,b]$.  It
  follows 
  that $d$ is a distance on $X$.

  Let $V \subset X$ be an open set (for the original topology) and let
  $x\in V $.
  Fix a chart $(U, \tau, A,\Omega)$ around $x$ and a local extension
  $\om'$.  Denote by $d_\Omega$ the Riemannian distance of $(\Omega,
  \om')$.  Let $U' $ be a neighbourhood of $x$ with compact closure in
  $U \cap V$. Since $\tau(x) \not \in \tau(\partial U')$,
  $\eps=d_\Omega\bigl (\tau(x), \tau(\partial U') \bigr) >0$.  If
  $\alfa: [a,b]\ra X$ is a continuous path with $\alfa(a) =x$ and
  $\alfa(b) \not \in U'$ set $c= \sup \{ t \in [a,b]: \alfa(t) \in
  U'\}$. Then
  \begin{gather*}
    L(\alfa) \geq L(\alfa \restr{[a,c]}) = L(\tau \circ \alfa
    \restr{[a,c]}) \geq d_\Omega \bigl (\tau(x), \tau(\partial
    U')\bigr ) = \eps.
  \end{gather*}
  Hence $\Bd(x, \eps) \subset U' \subset V$. This shows that the
  metric topology is finer than the original one.


  Conversely we show that for any $x \in X$ and $\delta >0$ the metric
  ball $\Bd(x, \delta)$ is open in the original topology.  Let again
  $(U, \tau, A,\Omega)$ be a chart around $x$ and let $\om'$ be a
  local extension and assume that there is a constant $ C>0$ such that
  $C\meno |d\tau( v) | \leq |v|_g \leq C |d\tau(v)| $ for any $v \in
  TU_\reg$.  Thanks to \L ojasiewicz Theorem by restricting $U$ and
  $\Omega$ we can assume that for any $z,z'\in A$ there is a $\cl^1$
  path joining $z$ and $z'$ and having Euclidean length $\leq C'
  |z-z'|^\mu$.  For $x' \in \Bd(x,\delta)$ put $\delta' = \sqrt[\mu]{(
    \delta - d(x,x'))/ CC'} >0$. Then the set $\tau\meno(\{ z\in
  \Omega: |z - \tau(x')| < \delta'\}) $ is contained in $ \Bd(x',
  \delta - d(x,x')) \subset \Bd(x, \delta)$. Therefore $\Bd(x,
  \delta)$ is open in the original topology and the two topologies
  coincide.
\end{proof}
Starting from the metric space $(X,\dom)$ one can define a new length
functional $L_d $ by the formula
\begin{equation}
  \label{eq:def-L-dom}
  L_{\dom}(\gamma) = \sup \sum_{i=1}^N \dom(\gamma(t_{i-1}) \gamma(t_i))
\end{equation}
the supremum being over all partitions $t_0 < \ldots{} < t_N$ of the
domain of $\gamma$.  By definition $d(x,y) \leq L_d(\ga) $ for any
continuous path joining $x$ to $y$, while the inequality $L_d(\ga)
\leq L(\ga)$ holds for any piecewise $\cl^1$ path.  The distance
$\dom$ is \enf{intrinsic} if $ \dom(x,y)= \inf \{L_\dom (\gamma):
\gamma \in C([0,1],X), \gamma(0)=x, \gamma (1)=y\}.  $
\begin{prop}
  The distance $d$ is intrinsic.
\end{prop}
\begin{proof}
  This is proved for general length structures in \cite[Prop. 2.4.1
  p.38]{burago-burago-ivanov}.
\end{proof}
\begin{definiz}
  We call $d$ the \enf{intrinsic distance} of the \Keler curve
  $(X,\om)$.
\end{definiz}
For geodesics in the metric space $(X,d)$ we adopt the following
terminology.  A \enf{shortest path} is a map $\ga : [a,b]\ra X $ such
that $L_d(\ga) = d(\ga(a), \ga(b))$.  \enf{Minimising geodesic} is
synonymous of shortest path.  One can reparametrise a shortest path in
such a way that $\dom(\ga(t), \ga(t'))=|t-t'|$ for any $t,t'$. In this
case we say that $\ga$ has \enf{unit speed}.
A \enf{segment} is by definition a unit speed shortest path.  More
generally, we say that $\ga$ is parametrised with constant speed $c$
if $\dom(\ga(t), \ga(t'))=|t-t'|$ for any $t,t'$.  If $I$ is any
interval a path $\ga: I\ra X$ is a \enf{geodesic} if for any $t \in I$
there is a compact neighbourhood $[t_0, t_1]$ of $t$ in $I$ such that
$\ga\restr{[t_0, t_1]}$ is a shortest path with constant speed.
\begin{lemma}[ \protect {\cite[Prop. 2.5.19
    p.49]{burago-burago-ivanov}}]
  \label{ball-hopf} If the ball $\Bd(x,r)$ is relatively compact in
  $X$, for any $y\in \Bd(x,r)$ there is a segment
  from $x$ to $y$.
\end{lemma}
$(\Xr, g)$ is a (smooth) Riemann surface with a \enf{noncomplete}
smooth \Keler metric.  For $x\in\Xr$ denote by $U_xX $ the unit sphere
in $T_xX$. Let $U\Xr=\bigcup_{x\in \Xr} U_xX$ be the unit tangent
bundle. We denote by $ (t, v) \mapsto \ga^v(t)$ the geodesic flow:
that is $\ga^v(t) = \exp_x(tv)$ where $x=\pi(v)$.  Let $\mathscr{U}
\subset \R\times T\Xr$ be the maximal domain of definition of the
geodesic flow of $(\Xr, g)$.  It is an open neighbourhood of
$\{0\}\times T\Xr$ in $\R\times T\Xr$.  Let $\Dom \subset T\Xr$ denote
the maximal domain of definition of the exponential: $\Dom = \{ v\in
T\Xr: (1,v) \in \mathscr{U}\}$.  For $x\in \Xr$ set $\Dom_x = \Dom
\cap T_x X$. Then $\mathscr{D}_x$ is the maximal domain of definition
of $\exp_x$. Both $\Dom$ and $\mathscr{D}_x$ are open in $T\Xr$ and
$T_xX$ respectively and the maps $\exp: \Dom \ra \Xr$ and $\exp_x :
\Dom_x \ra \Xr$ are defined and smooth.  For 
$v\in U_xX$ set
\begin{gather}
  \label{eq:def-Tv}
  T_v=\sup \{ t>0: tv \in \Dom_x\}.
\end{gather}
Denote by $\bx(0,r)$ the ball in $T_xX$ with respect to $g_x$.
\begin{definiz}
  For $x\in X_\reg$ the \enf{injectivity radius} at $x$, denoted
  $\inj{x}$, is the least upper bound of all $\delta>0$ such that
  $\bx(0,\delta) \subset \Dom_x$ and $ \exp_x\restr{\bx(0,\delta)} $
  is a diffeomorphism onto its image.
\end{definiz}
\begin{lemma}
  \label{inj-distanza-da-sing}
  For any $x\in \Xr$, $\inj{x} \leq \dom(x, X_\sing)$.
\end{lemma}
\begin{proof}
  Let $\ga : I \ra X$ be a piecewise $\cl^1$ path in $\A$ joining $x$
  to some singular point $x_0$. For $\delta \in (0, \inj{x})$ put
  $U_\delta = \exp_x(\bx(0,\delta))$.  Since $\overline{U_\delta}
  \subset X_\reg$ and $x_0$ is singular, there is some $t\in I$ such
  that $\gamma(t) \in \partial U_\delta$.  Let $t_0$ be the smallest
  such number. Then $\gamma\restr{[0,t_0]}$ is a path entirely
  contained in $X_\reg$. It follows from Gauss Lemma \cite[Prop. 3.6,
  p.70]{do-carmo} that $L(\gamma\restr{[0,t_0]}) \geq
  \delta$. Therefore also $L(\gamma)\geq \delta$. Since $\gamma$,
  $x_0$ and $\delta <\inj{x}$ are arbitrary we get $\dom(x, X_\sing)
  \geq \inj{x}$.
\end{proof}

\begin{lemma}\label{lemma-enumerato}
  Let $x\in \Xr $ and $y\in \Bd(x,\inj{x})$.
  \begin{enumerate}
  \item The intrinsic distance equals the Riemannian distance in
    $(\Xr, g)$:
    \begin{equation}
      \label{d=reg}
      \begin{aligned}
      d(x,y)=\inf\{L(\ga):
 &\      \ga \text { piecewise } \cl^1 \text { path in } \Xr  \\
&\text { with }
      \ga(0)=x, \ga(1)=y\}.
      \end{aligned}
    \end{equation}
  \item $\Bd(x,\inj{x}) = \exp_x( \bx(0,\inj{x}))$.
  \item There is a unique segment joining $x$ to $y$ and it coincides
    with the minimising Riemannian geodesic in $(\Xr, g)$ from $x$ to
    $y$.
  \item \label{poppa} A geodesic $\ga $ in $(X,d)$ is smooth on
    $\ga\meno (\Xr)$ and there $\nabla_{\dga} \dga =0$.
  \end{enumerate}
\end{lemma}
\begin{proof}
  It follows from the previous lemma and the hypothesis $d(x,y)<
  \inj{x}$ that paths passing through singular points do not
  contribute to the infimum in \eqref{eq:d-om}. This proves
  \eqref{d=reg}. From this follows that $y$ lies in
  $\exp_x(\bx(0,\inj{x}))$. So $ \Bd(x,\inj{x}) \subset
  \exp_x(\bx(0,\inj{x})) $. The reverse implication is obvious.  This
  proves 2.  In particular $y\in \exp_x(\bx(0,\inj{x}))$, so there are
  $v\in U_xX$ and $r\in (0, \inj{x})$ such that $y=\exp_x rv$. It
  follows from Gauss Lemma that the $\inf$ in \eqref{d=reg} is
  attained only on the path $\ga(t)= \exp_x t v$, $t\in [0,r]$. So
  $L(\ga)=d(x,y)$. But $ d(x,y) \leq L_d(\ga)\leq L(\ga)$ so $L_d(\ga
  ) = L(\ga)$ and $\ga$ is a segment also in $(X,d)$.  We have to
  prove that it is the unique one.  Since $d(x,y) < \inj{x}$ it
  follows from 2 that any other segment $\alfa$ must lie in
  $\exp_x(\bx(0,\inj{x})) \subset \Xr$. If $\alfa$ is smooth we can
  again apply Gauss Lemma. So it is enough to show that $\alfa$ is
  differentiable, which will yield 4 at once. This is a local problem,
  so we just prove that $\alfa\restr{[0, t_0]}$ is smooth for some
  $t_0>0$.  By Whitehead theorem \cite[Prop.  4.2 p.76]{do-carmo}
  there is a neighbourhood $W$ of $x$ such that for any $z\in W$,
  $W\subset \Bd(z,\inj{z})$. Let $t_0$ be small enough so that
  $\alfa([0,t_0]) \subset W$.  Put $x_0=\alfa(t_0)$ and let $\beta$ be
  the unique minimising Riemannian geodesic from $x$ to $x_0$. We
  already know that $L(\beta) = d(x,x_0) = t_0$.  For $t\in (0,t_0)$
  let $\beta_1$ and $\beta_2$ be the unique Riemannian geodesics from
  $x $ to $\alfa(t)$ and from $\alfa(t)$ to $x_0$ respectively. Both
  of them are also shortest paths, by the above.  Moreover
  \begin{equation*}
    t_0 = L_d(\alfa) \geq d(x,\alfa(t)) + d(\alfa(t), x_0) = L(\beta_1) +
    L(\beta_2) \geq L(\beta) = t_0. 
  \end{equation*}
  So $ L (\beta_1 * \beta_2)=L(\beta_1) + L(\beta_2)= L(\beta)$. Since
  the concatenation $\beta_1 *\beta_2$ is piecewise smooth $\beta_1 *
  \beta_2 = \beta$. This means that $\alfa(t)$ lies on
  $\beta([0,t_0])$. Since $t$ is arbitrary we get
  $\alfa\restr{[0,t_0]}=\beta$. In particular $\alfa$ is smooth.
\end{proof}

\begin{prop}
  On piecewise $\cl^1$ paths the functional $L_\dom$ agrees with $L$.
\end{prop}
\begin{proof}
  The inequality $L_d \leq L$ is obvious from the definition of $d$.
  For the reverse inequality consider a piecewise $\cl^1$ path $\gamma
  :[0,1]\ra X $ and assume at first that $\ga([0,1]) \subset \Xr$.
  Since $L_d(\ga) \leq L(\ga)<\infty$ the limit
  \begin{equation*}
    \lim_{h\ra 0} \frac{d(\ga(t+h),\ga(t))}{|h|}.
  \end{equation*}
  exists for a.e. $t\in [0,1]$. It is called \enf{metric derivative}
  and denoted by $ |\dga(t)|_d $. It is an integrable function of $t$
  and
  \begin{equation*}
    L_d(\ga) = \int_0^1 |\dga(t)|_d \ dt.
  \end{equation*} 
  (See \cite[p.106-109]{rinow} or \cite[p.59ff]{ambrosio-tilli}.)  So
  it is enough to check that $|\dga|_d = |\dga|_g$. This is
  accomplished as follows. Put $x=\ga(t)$. For small $h$ we can write
  $\ga(t+h) = \exp_x (z(h))$ where $z=z(h)$ is some $\cl^1$ path in
  $T_xX$ with $z(0)=0$ and and $ \dot{z}(0)= d(\exp_x)_0 (\dot{z}(0))=
  \dga(t)$. Since $d(\ga(t+h), \ga(t))= |z(h)|_g$
  \begin{equation*}
    |\dga(t)|_d = 
    \lim_{h\ra 0} \frac{|z(h)|_g}{|h|} = |\dot{z}(0)|_g= |\dga(t)|_g.
  \end{equation*}
  This proves that $L_d=L$ for paths that do not meet $X_\sing$.
  (There is a proof for $\cl^1$ Finsler manifolds due to Busemann and
  Mayer. It can be found in \cite{busemann-mayer} or at pp. 134-140 of
  Rinow's book \cite{rinow}.)  For a general path one can reason as in
  Lemma \ref{vero-length}: let $E=\gamma^\menuno{}(X_\reg)$,
  $F=I\setminus E$, $B=F^0$, $D=\partial F$. Then $I=E \sqcup B \sqcup
  D$. Since $\gamma$ is constant on the connected components of $F$,
  if $[a,b]$ is one such component then $L_d(\ga_{[a,b]} =
  L(\ga_{[a,b]})= 0 $.  The result follows from additivity of both
  functionals.  \end{proof}

\begin{cor}\label{semicont}
  The functional $L$ is lower semicontinuous on the set of piecewise
  $\cl^1$ paths with respect to the topology of pointwise convergence.
\end{cor}
\begin{proof}
  It easily follows from the definition that $L_d$ is lower
  semicontinuous on $\cl^0([0,1],X)$ with respect to pointwise
  convergence \cite[Prop. 2.3.4(iv)]{burago-burago-ivanov}.
\end{proof}

The construction of the intrinsic distance is \enf{local} in the
following sense.
\begin{lemma}\label{localize}
  For any point $x_0 \in X$ and any neighbourhood $U$ of $x_0$ in $X$
  there is a smaller neighbourhood $U'\subset U$ such that for any $x,
  y \in U'$ there is a segment $\gamma$ from $x$ to $y$ and any such
  segment is contained in $U'$.  In particular the intrinsic distance
  of $(X,\om)$ and that of $(U, \om\restr{U})$ coincide on $U'$.
\end{lemma}
\begin{proof} Let $\eps>0$ be such that $\overline{\Bd(x_0, 4\eps)}$
  is a compact subset of $U$.  Put $U'=\Bd(x,\eps)$. If $x,y \in U'$
  then $d(x, y)\leq d(x, x_0) + d(x_0, y) < 2\eps$.  By Lemma
  \ref{ball-hopf} since $\overline{\Bd(x, 2\eps)} $
  is compact there is a segment from $x$ to $y$.  Now if
  $\gamma=\gamma(t)$ is \enf{any} such segment $d(\gamma(t),x_0) \leq
  d(\gamma(t), x) + d(x, x_0) \leq L(\gamma) + d(x, x_0) \leq 3\eps$.
  So $\ga(t) $ lies in $U$.
\end{proof}

%
%

\begin{cor}\label{localizzo}
  Let $(X,\om)$ be a \Keler curve and let $d$ be the intrinsic
  distance.  If $(U, \tau, A,\Omega)$ is a chart around $x \in X$ and
  $\om'$ is a local extension of $\om$ there is a neighbourhood $U'
  \subset U$ of $x$ such that $\tau\restr{U'}$ is a biholomorphic
  isometry between $(U',d)$ and $\tau(U') \subset A$ provided with the
  intrinsic distance obtained from $\om'$.
\end{cor}

It follows that to study \enf{local} properties of the metric spaces
$(X,d)$ it is enough to consider the special case in which $X$ is an
analytic set in a domain of $\C^n$ with the metric induced from some
\Keler metric of the domain.  This situation, under the additional
hypothesis that the singularity be analytically irreducible, is the
object of \Sec\Sec \ref{sec-regu}--\ref{convex-section}, throughout
which we will make the following assumptions and use the following
notation.
\begin{quote}
  \label{notazione}
  $\langle \, \, , \, \rangle$ is the standard Hermitian product on
  $\C^n$, \\
  $v\cdot w = \Re \langle v,
  w\rangle $ is the corresponding scalar product,\\
  $|\cdot|$ is the
  corresponding norm.\\
  Given two nonzero vectors $v, w$ in a Euclidean space
  \begin{equation*}
    \ango(v,w)=\arccos\frac{v\cdot w}{|v|\cdot |w|}
  \end{equation*}
  is the unoriented angle between them.\\
  $\Omega'$ is an open
  polydisc centred at $0\in \C^n$,\\
  $A\subset \Omega'$ is an analytic
  curve,\\
  $\om$ is a smooth \Keler form on $\Omega'$,
  \\
  $g$ is the corresponding \Keler metric,
  \\
  $g_x$ is the value of $g$ at $x \in \Omega'$,\\
  $|\cdot|_g$ or  $|\cdot|_x$ denotes the corresponding norm,\\
  $g_0=\langle \, , \, \rangle$, \\
  $\Omega$ is an open subset of $\Omega'$ with
  $\overline{\Omega} \compsubset{} \Omega'$,\\
  $X:=A\cap \Omega$,\\
  $d$ is the intrinsic distance of $(X,\om\restr{X})$,\\
  $\Bd(x,r)$ is the  ball in $(X,d)$,\\
  $\Bd^*(x,r) = \Bd(x,r) \setminus \{0\}$,\\
  $B_x(0,r) = \{ w\in T_xX : |w|_x< r\}$.\\
  $X_\sing = \{0\}$, \\
  $X$ is analytically irreducible
  at $0$,\\
  $m=\mult_0 X$ is the multiplicity of $X$ at $0$,\\
  $K(x)$ is the Gaussian curvature of $(\Xr, g)$ at $x\in \Xr$ and\\
  \begin{equation}
    \label{eq:def-K0}
    \Kl = \sup_{x\in\Xr} K(x).
  \end{equation}
  $\Delta=\{z\in \C: |z| < 1\}$,\\
  $\Delta^*=\Delta\setminus\{0\}$,\\
  $\Delta' \subset \C$ is an open subset containing $\overline{\Delta}$,\\
  $\phi : \Delta' \ra X'$ is the normalisation map,\\
  $\phi(\Delta) = X$. \\
  There is a holomorphic map $\psi=(\psi_1, \ldots{}, \psi_n): \Delta'
  \ra \C^n$ such that
  \begin{equation}
    \label{normalizzazione}
    \phi(\z) = \z^m \psi(\z) \qquad
    \psi_1(z)\equiv 1 \qquad \psi_j(0)=0 \  j>1.
  \end{equation}
  $R: \Delta' \ra \C^n$ is the holomorphic map defined by
  \begin{gather}
    \label{eq:def-Rz}
    R(\z) := \frac{\psi(\z) - \psi(0)}{m \z} + \frac{\psi'(\z)}{m}
    \\ \label{eq:phi-primo-con-Rz} \phi'(\z) = m \z^{m-1} (e_1 + \z
    R(\z))
  \end{gather}
  $e_1=(1,0, \ldots{}, 0)$.\\
  $c_0>0$ is a constant such that
  \begin{gather}
    \label{eq:def-c0}
    \sup_\Delta |\phi'| \leq c_0 \qquad \sup_{\Delta}|R|\leq c_0 \\
    \label{eq:c0-metrico}
    \forall x\in \Omega, \forall v\in \C^n, \quad
    \begin{cases}
      \tfrac{1}{c_0} |v| \leq |v|_x \leq c_0 |v| \\
      |v|_{x} \leq |v| ( 1 + c_0 |x|).
    \end{cases}
  \end{gather}
  From \eqref{eq:def-c0} it follows that for any $z\in \Delta$
  \begin{gather}
    \label{eq:c0-zeta}
    |\phi(z)| \leq c_0 |z|.
  \end{gather}
  $\pi: \C^n \ra \C \times \{0\}$ is the projection on the first coordinate,\\
  $ u:= \pi \circ \phi : \Delta \ra \Delta $ is the standard $m:1$
  ramified covering: $u(z)=z^m$.\\
  For $\theta_0 \in \R$ and $\alpha \in(0,\pi]$ put
  \begin{equation}
    \label{eq:def-sectors}
    S(\theta_0, \alpha) = \{ \rho e^{i \theta} : \rho\in (0,1), |\theta
    - \theta_0|< \alpha\} \subset \Delta.
  \end{equation}
  Then
  \begin{equation}
    \label{eq:controsettori}
    u^\men \bigl ( S(\theta_0, \alpha)\bigr ) = \bigsqcup_{j=0}^{m-1}
    S\biggl (\frac{\theta_0 + 2\pi j}{m} , \frac{\alpha}{m} \biggr )
  \end{equation}
  and
  \begin{equation}\label{eq:def-u_j}
    u_j:= u \Big \vert  _{  S\bigl (\tfrac{\theta_0 + 2\pi j}{m} ,
      \tfrac{\alpha}{m} \bigr ) }
  \end{equation}
  is a biholomorphism onto $S(\theta_0, \alpha)$.  The Whitney tangent
  cone of $X$ at $0$ is
  \begin{gather}
    \label{eq:Whitney}
    C_0 X =\C \times \{0\} \subset \C^n
  \end{gather}
  (see e.g.  \cite[p.122, p.80]{chirka}). \\
  If $\ga : [0,L] \ra X$ is a path,
  $\ga\no (t) = \ga(L-t)$.
\end{quote}

\section {Regularity of geodesics}
\label{sec-regu}

\begin{lemma} [\protect{\cite[Lemma 1, p.86]{lojasiewicz-ihes}}]
  \label{loja-radici}
  Let $m$ be a positive integer and $K>0$.  Put $ Z=\{(a_1,
  \ldots{},a_m,x) \in \C^{m+1} : \ x^m + \sum_{j=1}^m a_j x^{m-j} =0,\
  |a_j|\leq K\}$.  Then there is an $M=M(m,K)>0 $ with the following
  property. Let $\alpha(t) = (a(t), x(t))$ be a continuous path
  $\alpha: [0,1] \ra Z$ and $L>0$ such that $ | a(t) - a(t')| \leq L
  |t-t'| $ for $t,t' \in [0,1]$.  Then
  \begin{equation}
    |x(t) - x(t')| \leq M L^{\nicefrac{1}{m}}  |t-t'| ^\um
    \qquad \forall t, t'\in [0,1].
  \end{equation}
\end{lemma}

\begin{prop}\label{param-hoelder}
  There is a constant $\cuno> 1$ such that for any $z,z'\in \Delta$
  \begin{equation}
    \label{eq:hold-phi}
    \frac{1}{\cuno} d(\phi(z), \phi(z')) \leq
    |z-z'| \leq \cuno d(\phi(z), \phi(z'))^\um .
  \end{equation}
\end{prop}
\begin{proof}
  Recall that $\Omega'$ is a polydisc, say $\Omega' =P(0)_{K,\ldots{},
    K}$ and $X$ is compactly contained in $\Omega'$.  Let $z, z'\in
  \Delta$ and $x=\phi(z), x'=\phi(z')$. For $\eps>0$ let $\ga: [0,1]
  \ra X$ be a piecewise $\cl^1$ path with $L:=L(\ga)< d(x,x') + \eps$.
  We can assume that $\ga$ has constant speed equal to $L$, so
  $d(\ga(t), \ga(t')) \leq L |t-t'|$.  On the other hand we trivially
  have $|\ga(t) - \ga(t')| \leq d(\ga(t), \ga(t'))$.  Put $a_m(t) =
  -\pi\bigl ( \ga(t) \bigr )$, $x(t) = \phi\meno(\ga(t))$ and
  $\alpha(t) = (0,\ldots{},0, a_m(t),x(t))$.  Then
  \begin{gather*}
    a_m(t)=- \pi\circ \phi (x(t))= - u(x(t)) = - x^m(t) \qquad x^m(t)
    +
    a_m(t) =0 \\
%
    |a_m(t) - a_m(t')| = |\pi(\ga(t)) - \pi(\ga(t'))| \leq \\
    \leq |\ga(t) - \ga(t')| \leq d (\ga(t) , \ga(t')) \leq L|t-t'|.
  \end{gather*}
  Therefore by Lemma \ref{loja-radici} applied to $\alfa$
  \begin{equation*}
    |x(t)- x(t')| \leq M L^\um |t-t'|^\um.
  \end{equation*}
  For $t=0$ and $t=1$ we get
  \begin{equation*}
    |z-z' | \leq  M L^\um \leq M  \bigl (
    d(x, x') + \eps \bigr )^\um.
  \end{equation*}
  Letting $\eps \ra 0$ we get $ |z-z'| \leq M d(\phi(z), \phi(z'))^\um
  $.  On the other hand it follows from
  \eqref{eq:def-c0} that $d(\phi(z), \phi(z') ) \leq c_0 |z-z'|$, so
  $\cuno= \max \{c_0, M\}$ works.
\end{proof}

\begin{cor}
  For any $r$ with $ 0 < r< \cuno^{-m}$
  \begin{gather}
    \label{eq:palletta-dischetto}
    \Bd(0,r) \subset \phi(B(0, \cuno r^\um)) \subset \Bd(0, \cuno^2
    r^\um ).
  \end{gather}
\end{cor}

For $x\in \Xr$ let $(T_xX)^\perp$ denote the orthogonal complement of
$T_xX \subset \C^n$ with respect to the scalar product $g_x$. If $w\in
\C^n$, $w^\perp$ denotes the $g_x$--orthogonal projection of $w$ on
$(T_xX)^\perp$.  Let $\sef_x : T_xX \times T_xX \ra (T_xX)^\perp$ be
the second fundamental form of $\Xr$. Since $g$ is \Keler and $X_\reg$
is a complex submanifold $\sef_x$ is complex linear.  If $v$ is a
nonzero vector in $T_xX$ put
\begin{equation}
  |\sef_x| = \frac{|\sef_x(v,v)|_x }{|v|_x^2}.
\end{equation}
Since $T_xX$ is complex one-dimensional the choice of $v$ is
immaterial.  Denoting by $K_\Omega(T_xX)$ the sectional curvature of
$(\Omega, g)$ on the 2-plane $T_xX$, Gauss equation yields
\begin{equation}
  \label{eq:Gauss}
  K(x)= K_\Omega(T_xX) - 2 |\sef_x|^2.
\end{equation}
(See e.g. \cite{KN-II} p. 175-176.)
\begin{prop}
  There is a constant $\cdue$ such that
  \begin{gather}
    \label{eq:II-est-z}
    |\sef_{\phi(z)}| \leq \frac{\cdue}{|z|^{m-1}} \qquad \forall z \in
    \Delta
    \\
    \label{eq:II-est-d}
    |\sef_{x}| \leq \frac{\cdue}{ d(x, 0) ^{1-\um}} \qquad \forall
    x\in X_\reg.
  \end{gather}
\end{prop}
\begin{proof}
  By \eqref{normalizzazione} we have $\phi(z)=z^m \psi(z)$, so
  $\phi'(z)= z^{m-1} v(z)$ where $v(z)= m\psi(z) + z \psi'(z)$.  Since
  $\phi$ is a holomorphic immersion on $\Delta'\setminus \{0\} $,
  $v(z) \neq 0 $ and $T_{\phi(z)} X = \C \phi'(z) = \C v(z)$ for any
  $z\neq 0$. But also $v(0)=m e_1\neq 0$.  So $v : \Delta' \ra \C^m$
  is continuous and nonvanishing,
  hence $ \inf _\Delta |v|_g >0 $ and $ \sup _\Delta |v|_g < +\infty$.
  Similarly $ \sup_\Delta |v'|_g < +\infty$.  Let $C$ be such that
  \begin{gather*}
    \inf _\Delta |v|_g \geq \frac{1}{C} \qquad \sup _\Delta |v|_g \leq
    C \qquad \sup_\Delta |v'|_g \leq C.
  \end{gather*}
  Then we have
  \begin{gather}
    \nonumber{} \sef_{\phi(z)} (v(z),v(z)) = \frac{1}{z^{m-1}}
    \sef_{\phi(z)} (\phi'(z), v(z) )
    = \frac{1}{z^{m-1}} \bigl(v'(z)\bigr)^\perp\\
    \big|\sef_{\phi(z)} \big| = \frac{
      |(v'(z))^\perp|_{\phi(z)}}{|z|^{m-1} |v|_{\phi(z)}^2} \leq
    \frac{ |v'(z))|_{\phi(z)}}{|z|^{m-1} |v|_{\phi(z)}^2} \leq \frac{
      C^3}{|z|^{m-1}}.
    \label{eq:B-uguale}
  \end{gather}
  This proves \eqref{eq:II-est-z}. 
  For $x\in \Xr$ let $z=\phi^\men (x)$ and consider the path $\ga(t) =
  \phi(tz)$, $t\in [0,1]$. Then
  \begin{gather*}
    \dga(t) = \phi'(tz)z = (tz)^{m-1} v(tz) z = z^m t^{m-1} v(tz)\\
    |\dga(t)|_x = |z|^m t^{m-1} |v(tz)|_x \leq C |z|^m\\
    d(x,0) \leq \Lo(\ga) = \int_0^1 |\dga(t)|\ dt \leq C |z|^m\\
    |z|\geq 
    \sqrt[m]{\frac{ \dom(x,0)}{C}}.
  \end{gather*}
  So \eqref{eq:II-est-d} follows from \eqref{eq:II-est-z}.
\end{proof}

\begin{remark} The map $v$ above is a holomorphic vector field along
  the map $\phi: \Delta \ra X$. On the other hand the push forward of
  $v$ to $X$, that is $v\circ \phi^\men$, is only weakly
  holomorphic. In fact any holomorphic vector field on $X$ has to
  vanish at $0$ if $X$ is singular
  \cite[Thm. 3.2]{rossi-vector-fields}.
\end{remark}

\begin{lemma}\label{anti-Hol}
  If $a, b \geq 0$ and $s\in (0,1)$ then $ |a^s-b^s|\leq |a-b|^s$.
\end{lemma}
\begin{proof} Assume $a\geq b$.  The function $\eta(x) = (b+x)^s -
  x^s$ belongs to $ \cl^0([0, +\infty))\cap \cl^1((0,+\infty))$.
  Since $s<1$, $\eta'(x) = s[ (b+x)^{s-1} -x^{s-1}] \leq 0$. So
  $\eta(a-b) = a^s - (a-b)^s \leq b^s$.  \end{proof}

\begin{teo}\label{W-Hoelder}
  There is a constant $\ctre$ such that for any unit speed geodesic
  $\ga : [0,L]\ra X$ with $\ga((0,L]) \subset X_\reg$ we have
  \begin{equation}
    \label{eq:hoelder}
    ||\dga||_{\cl^{0, \nicefrac{1}{m}}}  \leq \ctre
  \end{equation}
  Here the H\"older norm is computed using the Euclidean distance on
  $\C^n$.
\end{teo}
\begin{proof}
  By \ref{poppa} of Lemma \ref{lemma-enumerato} $\ga\restr{(0, L]}$ is
  a Riemannian geodesic of $X_\reg$. Hence the acceleration $\ddga(t)$
  is orthogonal to $T_{\ga(t)}X$ with respect to the scalar product
  $g$. So for $t>0$, $ \ddga(t) = \bigl ( \ddga (t)\bigr )^\perp =
  \sef_{ \ga(t)} \bigl ( \dga(t), \dga(t)\bigr ) $.  Using
  \eqref{eq:II-est-d}, $|\dga|\equiv 1$ and \eqref{eq:c0-metrico}
  we get
  \begin{equation*}
    \begin{gathered}
      |\ddga(t)|\leq \czero |\ddga(t) |_x = \czero \big |\sef_{
        \ga(t)} \bigl ( \dga(t), \dga(t)\bigr ) \big |_x =\czero \big
      |\sef_{\ga(t)} \big| \leq \frac{C} {d(\ga(t), 0)^{1-\um}}
    \end{gathered}
  \end{equation*}
  where $C=\czero{} \cdue$.  Set $a:=d(\ga(0), 0)$ and
  $\beta=1-1/m$. For $t>0$ we have
  \begin{gather}
    \nonumber{} d(\ga(t), 0) \geq \big |d(\ga(t), \ga(0)) - d(\ga(0),
    0) \big| = |t-a|
    \\
%
%
    \label{eq:stima-ddga}
    |\ddga(t)|\leq \frac{C} {|t-a|^{\beta}}.
  \end{gather}
  We claim that
  \begin{equation}
    \label{eq:stima-conclamata}
    |\dga(t) - \dga(s) | \leq 2mC |t-s|^{\um}
  \end{equation}
  for any pair of numbers $s, t$ such that $0<s\leq t\leq 1$. Indeed
  \begin{equation*}
    |\dga(t) - \dga(s) | \leq  
    \int_s^t  |\ddga(\tau)|d\tau \leq 
    C  \int_s^t \frac{d\tau}{ |\tau-a|^\beta} = C \bigl  ( I_a(t) -
    I_a(s) \bigr ) 
  \end{equation*}
  where we put $I_a(t) = \int_0^t |\tau-a|^{-\beta}d\tau$.  A simple
  computation shows that
  \begin{equation*}
    I_a(t) -
    I_a(s) =
    \begin{cases}
      m \bigl ( |s-a|^{\um} - |t-a|^{\um} \bigr ) & \text{ for } 0<s
      \leq t   \leq a \\
      m \bigl ( |s-a|^{\um} + |t-a|^{\um} \bigr ) & \text{ for } 0<s
      \leq a   \leq t \\
      m \bigl ( |t-a|^{\um} - |s-a|^{\um} \bigr ) & \text{ for } a < s
      \leq t    \\
    \end{cases}
  \end{equation*}
  By Lemma \ref{anti-Hol}
  \begin{gather*}
    \Bigl ||t-a|^\um - |s-a|^{\um} \Bigr | \leq \Bigl | |t-a| - |s-a|
    \Bigr |^{\um} \leq |t-s| ^\um
  \end{gather*}
  so $I_a(t) - I_a(s) \leq m |t-s|^{\um}$ in the first and the last
  case. As for the middle case, namely $s\leq a \leq t$, we have
  $|s-a|\leq |s-t| $ and $|t-a|\leq |t-s|$, so $ |s-a|^{\um} +
  |t-a|^{\um} \leq 2|t-s|^{\um}$.  Therefore in any case $ I_a(t) -
  I_a(s) \leq 2m |t-s|^{\um} $ and this finally yields
  \eqref{eq:stima-conclamata}.  This proves \eqref{eq:hoelder} with
  $\ctre=2mC=2m\czero\cdue$.
\end{proof}

\begin{cor}\label{vecteur-tangent}
  Let $\ga : [0,L] \ra X$ be a segment with $\ga(0)=0$. Then $\ga$ is
  differentiable at $t=0$ and the map $\dga:[0,L] \ra \C^n$ is a
  H\"{o}lder continuous of exponent $\nicefrac{1}{m}$.
\end{cor}
\begin{proof}
  Since shortest paths are injective $\ga((0,L])\subset \Xr$. So
  estimate \eqref{eq:hoelder} holds.  Therefore $\dga$ is uniformly
  continuous on $(0,L]$ and extends continuously for $t=0$.  By the
  mean value theorem the extension for $t=0$ is precisely the
  derivative $\dga(0)$.
\end{proof}

%

\begin{lemma} \label{angoletto-1} For any $\eps>0$ there is a
  $\delta>0$ such that for any $x\in \Bd^*(0, \delta) $ and any $v\in
  T_xX$
  \begin{gather}
    \label{eq:sottosopra}
    (1-\eps)|v| < |v|_x < (1+\eps) |v|
    \\
    |\pi (v) - v| < \eps |v|\\
    (1-\eps) |v| < |\pi(v)| < (1+\eps) |v|.
  \end{gather}
\end{lemma}
\begin{proof}
  \eqref{eq:sottosopra} holds for $x$ sufficiently close to $0$ simply
  because $g_0=\langle\ , \, \rangle$.  For the second condition set
  \begin{gather*}
    \delta = 
    \eps^m\Bigl [{ (\cuno(1+c_0)(1+\eps)} \Bigr]^{-m}
  \end{gather*}
  where $c_0$ is the constant defined in \eqref{eq:def-c0}.  By Lemma
  \ref{param-hoelder} if $x\in \Bd(0,\delta)$ and $z =\phi\meno(x)\in
  \Delta$
  \begin{equation*}
    |z|< \cuno \delta^\um = \frac{\eps}{(1+c_0) (1+\eps)}.
  \end{equation*}
  Therefore
  $$
  |zR(z)| <\frac{\eps}{1+\eps} .
  $$
  ($R$ is defined in \eqref{eq:def-Rz}.)  It follows from
  \eqref{eq:phi-primo-con-Rz} that $T _x X = \C \cdot \phi'(z) = \C
  \cdot (e_1 + zR(z))$, so any $v\in T_xX$ is of the form $v= \la(e_1
  + zR(z)) $ for some $\la \in \C$. Then
  \begin{gather*}
    |v| \geq |\la| - |\la z R(z)| \geq |\la| -
    \frac{|\eps\la|}{1+\eps} = \frac{1}{1+\eps}|\la| \qquad |\la|\leq
    (1+\eps) |v|
    \\
    |v-\pi(v)| =\min_{w\in \C\times \{0\}} |v-w| \leq |v-\la e_1 |
    = |\la z R(z)|< |\la| \frac{\eps}{1+\eps} \leq \eps|v|\\
    \bigg | |\pi(v)| - |v| \bigg | \leq |\pi(v) -v| < \eps |v|.
  \end{gather*}
\end{proof}

\begin{lemma}
  \label{distanze-asintotiche}
  We have
  \begin{gather}
    \label{eq:distanze-asintotiche}
    \underset{\parbox{3em}{\centering\scriptsize $x,y \ra 0$ \\$ x,y
        \in X $ }} {\liminf} \ \frac{d(x,y) }
    {|x-y|} \geq 1.
  \end{gather}
\end{lemma}
\begin{proof}
  Given $\eps>0$ let $\delta>0$ be such that \eqref{eq:sottosopra}
  holds for any $x\in \Bd^*(0, \delta) $ and any $v\in T_xX$. If
  $x,y\in \Bd^*(0,\delta/3)$ and $\alfa: [0,L]\ra X$ is a segment,
  then $\alfa([0,L]) \subset \Bd(0,\delta)$ and the set $J=\{t\in
  [0,L]: \alfa(t) = 0\}$ contains at most one point. For $t\not \in J$
  $ |\dalf(t)|_{\alfa(t)} 
  \geq (1-\eps)|\dalf(t)|.  $ Integrating on $[0,L]\setminus J$ yields
  \begin{gather*}
    d(x,y) = L(\alfa) \geq (1-\eps) \int_0^L |\dalf(t)|dt \geq (1-\eps) |x-y|\\
    \frac{d(x,y) }{|x-y|} \geq 1-\eps.
  \end{gather*}
\end{proof}
\begin{lemma}\label{angoletto-2}
  For any $\eps >0$ there is a $\delta>0$ such that for any $x\in
  \Bd^*(0, \delta)$ and any pair of nonzero vectors $v,w\in T_xX$
  \begin{equation*}
    |    \ango(\pi(v) , \pi(w) ) - \ango(v,w)|< \eps.
  \end{equation*}
  (The angle is computed with respect to $\langle \ , \, \rangle$.)
\end{lemma}
\begin{proof}
  The angle function $\ango : S^{2m-1}\times S^{2m-1} \ra \R$ is the
  Riemannian distance for the standard metric on unit the sphere. In
  particular it is Lipschitz continuous with respect to the Euclidean
  distance. So one can find $\eps_1>0$ with the property that that
  \begin{gather}
    \label{eq:angolo-Lip}
    |u_1-u_2| < \eps_1, |w_1-w_2|< \eps_1 \Rightarrow |\ango(u_1, w_1)
    - \ango(u_2, w_2) | < \eps.
  \end{gather}
  We can assume $\eps_1 <1$.  Choose $\delta>0$ such that for $x\in
  \Bd^*(0, \delta)$ and $v\in T_xX$
  \begin{gather*}
    |\pi(v) - v| < \frac{\eps_1}{2}|v|.
  \end{gather*}
  This is possible by Lemma \ref{angoletto-1}. Moreover $v\neq 0
  \Rightarrow\pi(v)\neq 0$, because $\eps_1 <1$.  Then for two nonzero
  vectors $v, w \in T_xX$, $x\in \Bd^*(0,\delta)$
  \begin{gather*}
    \bigg |\frac{v}{|v|} - \frac{\pi(v)}{|\pi(v)|} \bigg | \leq
    \frac{2 |v-\pi(v)|} {|v|} < \eps_1 \qquad \bigg |\frac{w}{|w|} -
    \frac{\pi(w)}{|\pi(w)|} \bigg | < \eps_1.
  \end{gather*}
  Together with \eqref{eq:angolo-Lip} this yields the result.
\end{proof}

\begin{lemma} \label{angoletto-3} For any $\eps>0$ there is a
  $\delta>0$ such that for any segment $\ga: [0,L]\ra \Bd(0,\delta)$
  with $\ga((0,L)) \subset \Xr$ and any $s,s'\in [0,L]$
  \begin{gather}
    |\dga(s)- \dga(s')| < \eps\qquad \ango(\dga(s),\dga (s')) < \eps .
  \end{gather}
  (The angle is computed with respect to $\langle \ , \, \rangle$.)
\end{lemma}
\begin{proof}
  Let $\eps_1>0$ be such that
  \begin{equation}
    \label{eq:eps1-novo}
    |u-w| <
    \eps_1 \Rightarrow \ango(u,w) < \eps \qquad \forall
    u,w\in S^{2m-1} .
  \end{equation}
  Choose $\delta>0$ such that
  $$
  \sqrt[m]{2 \delta} < \min \biggl \{ \frac{\eps_1}{ 2 c_0\ctre},
  \frac{\eps}{\ctre} \biggr\} .
  $$
  If $\ga: [0,L]\ra \Bd(0,\delta)$ is a segment with $\ga((0,L))
  \subset \Xr$ at most one of the points $\ga(0)$ and $ \ga(L)$
  coincides with the origin.  So the H\"older estimate
  \eqref{eq:hoelder} holds for $\ga$.  Then
  \begin{gather*}
    |\dga(s) - \dga(s')| \leq \ctre \sqrt[m]{L} \leq \ctre \sqrt[m]{2
      \delta} < \eps
    . 
  \end{gather*}
  This proves the first inequality.  From \eqref{eq:c0-metrico} it
  follows that
  \begin{gather*}
    \frac{1}{|\dga(s)|} \leq c_0
  \end{gather*}
  so
  \begin{gather*}
    \bigg |\frac{\dga(s)}{|\dga(s)|} - \frac{\dga(s')}{|\dga(s')|}
    \bigg | \leq \frac{2 |\dga(s)-\dga(s')|} {|\dga(s)|} \leq 2 \czero
    \ctre \sqrt[m]{2\delta} < \eps_1.
  \end{gather*}
  Coupled with \eqref{eq:eps1-novo} this yields the second inequality.
\end{proof}

\begin{lemma} \label{angoletto-4} For any $\eps>0$ there is a
  $\delta>0$ such that for any segment $\ga: [0,L]\ra \Bd(0,\delta)$
  with $\ga((0,L)) \subset \Xr$ and any $s,s'\in [0,L]$
  \begin{gather*}
    \ango(\pi(\dga(s)),\pi(\dga (s'))) < \eps.
  \end{gather*}
  (The angle is computed with respect to $\langle \ , \, \rangle$.)
\end{lemma}
\begin{proof}
  Let $\eps_1>0$ be such that
  \begin{equation}
    \label{eq:primo-quartino}
    |u-w| <
    \eps_1 \Rightarrow \ango(u,w) <\frac{ \eps}{3} \qquad \forall
    u,w\in S^{2m-1} .
  \end{equation}
  Let $\delta_1>0$ be such that for any $x\in \Bd^*(0, \delta_1)$ and
  any $v\in T_xX$
  \begin{gather}
    \label{eq:primo-terzo}
    |\pi(v) -v| < \frac{\eps_1|v|}{\czero}.
  \end{gather}
  Such a $\delta_1$ exists by Lemma \ref{angoletto-1}.  Next, by Lemma
  \ref{angoletto-3}, there is $\delta_2>0$ such that for any segment
  $\ga: [0,L]\ra \Bd(0,\delta_2)$ with $\ga((0,L)) \subset \Xr$ and
  any $s,s'\in [0,L]$
  \begin{gather}
    \label{eq:secondo-terzo}
    \ango(\dga(s),\dga (s')) < \frac{\eps}{3} .
  \end{gather}
  Set $\delta= \min\{\delta_1, \delta_2 \}$. 
  If $\ga: [0,L]\ra \Bd(0,\delta)$ is a segment with $\ga((0,L))
  \subset \Xr$ and $s\in [0,L]$, then by \eqref{eq:primo-terzo} and
  \eqref{eq:c0-metrico}
  \begin{gather*}
    \big |\pi(\dga(s))-\dga (s)\big | < \frac{\eps_1 |\dga(s)|}
    {\czero} \leq \eps_1 |\dga(s)|_{\ga(s)} = \eps_1
  \end{gather*}
  so by \eqref{eq:primo-quartino}
  \begin{gather*}
    \ango\bigl (\pi(\dga(s)),\dga (s)\bigr) < \frac{\eps}{3}.
  \end{gather*}
  Then using \eqref{eq:secondo-terzo} we get for arbitrary $s,s'\in
  [0,L]$
  \begin{gather*}
    \ango(\pi(\dga(s)), \pi(\dga(s'))) \leq\\
    \leq \ango\bigl (\pi(\dga(s)),\dga (s)\bigr) +
    \ango(\dga(s),\dga(s')) + \ango\bigl ( \dga(s'),\pi(\dga(s'))
    \bigr ) < \eps
  \end{gather*}
  as claimed.
\end{proof}

\section{Tangent vectors in the normalisation}
\label{regu-delta}

%
In this section we study the regularity properties of the preimage in
$\Delta$ of segments in $X$.
\begin{lemma}
  If $\ga: [0,L] \ra X$ is a segment
  with $\ga(0)=0$, the path $\phi\meno\circ \ga : [0,L] \ra \Delta$
  has finite length.
\end{lemma}
\begin{proof} We know from Cor. \ref{vecteur-tangent} that $\dga(0)$
  exists.  By \eqref{eq:Whitney} $\dga(0)=(\dga_1(0), 0,\dotsc , 0)$
  and $\dga_1(0)=e^{i\theta_0}$ for some $\theta_0\in [0,2\pi)$.
  There is $\eps>0$ such that $\ga_1((0,\eps])$ is contained
  in the sector $ S(\theta_0, \pi)\subset \Delta$ defined in
  \eqref{eq:def-sectors}.  We can write $ \ga_1(t) = \rho(t)
  e^{i\theta(t)} $ for appropriate functions $\rho, \theta \in
  \cl^1((0,\eps])$.  Put $\beta=\phi\meno\circ \ga$.
  $\beta((0,t])$ is contained in one of the connected components of
  $\phi\meno S(\theta_0, \pi)$ hence by \eqref{eq:controsettori} there
  is an integer $k$, $0\leq k \leq m-1$, 
  such that
  \begin{gather}
    \beta(t) = u_k\meno (\ga_1(t))=\rho^\um (t) e^{i\theta(t) /m}
    \xi_k
  \end{gather}
  where $\xi_k = e^{2\pi k/m}$.  Then
  \begin{gather}
    \nonumber{} \dot{\beta} = \frac{1}{m} \rho^{\um -1} ( \rho' + i
    \theta'\rho)
    e^{i \theta /m} \xi_k\\
    \label{eq:norma-beta}
    | \dot{\beta}| = \frac{1}{m} \rho^{\um -1} \sqrt{ (\rho')^2 + i
      (\theta' \rho)^2}
    \\
    \nonumber{} \lim_{t\ra 0} \rho(t) =
    \lim_{t\ra 0} |\ga_1(t)|=0\\
    \lim_{t\ra 0} \frac{\rho(t)}{t} = \lim_{t\ra 0} \bigg \vert
    \frac{\ga_1(t)}{t} \bigg \vert = |\dga_1(0)| =1.
    \label{eq:rho-primo}
  \end{gather}
  If we put $\rho(0)=0$ and $ \rho'(0)=1$, then $\rho\in
  \cl^1([0,\eps])$.  Also
  \begin{gather*}
    e^{i\theta_0} = \dga_1(0) = \lim_{t\ra 0} \frac{\ga_1(t)}{t}=
    \lim_{t\ra 0} \frac{\ga_1(t)} {\rho(t)} = \lim_{t\ra 0}
    e^{i\theta(t)}\\
    \lim_{t\ra 0} \theta(t) = \theta_0 + 2 N \pi \qquad N\in \Zeta.
  \end{gather*}
  Change $\theta$ by subtracting $2N\pi$ to it and put $\theta(0)
  =\theta_0$. Then $\theta\in \cl^0([0,\eps])$.  
  \begin{gather*}
    \dga_1 = \bigl ( \rho' + i \rho\theta' \bigr) e^{i\theta} \qquad
    \lim_{t\ra 0} \dga_1(t) = \dga_1(0) = e^{i\theta(0)}
    \\
    \Longrightarrow \quad \lim_{t\ra 0} \rho(t) \theta' (t)=0.
  \end{gather*}
  Since $\rho'(0)=1$ and $\rho \theta ' \ra 0$, we get from
  \eqref{eq:norma-beta} that $|\dot{\beta}| \leq C_1 \rho^{\um -1}
  \leq C_2 t^{\um -1}$. Therefore $ L=\int_0^\eps |\dot{\beta}| < +
  \infty $ and $\beta$ has finite length.
\end{proof}

\begin{definiz}
  \label{def-unga} If $\ga: [0,L] \ra X$ is a unit speed geodesic with
  $\ga(0)=0$ denote by $\unga:[0, \unl]\ra \Delta$ the arc-length
  reparametrisation of the path $\phi\meno\circ \ga : [0,L] \ra
  \Delta$ (with respect to the Euclidean metric on $\Delta$).
\end{definiz}

\begin{prop} \label{rego-in-D} If $\ga: [0,L] \ra X$ is a unit speed
  geodesic with $\ga(0)=0$, then $\unga \in \cl^1([0,\unl])$ and
  $\dga(0)=(\dunga(0)^m, \ldots, 0)$.
\end{prop}
\begin{proof}
  Put $\beta=\phi\meno\circ \ga$.  By the previous Lemma $\beta$ has
  finite length.  If we set $h(t) = \int_0^t |\dot{\beta}(\tau)| d
  \tau$ then $\unga(s) = \beta(h^\men(s))$. It is clear that $\unga
  \in \cl^0([0,\unl]) \cap \cl^1((0, \unl])$, but we have to check
  that $\unga$ is continuously differentiable at $s=0$.  This is not
  immediate since $h'(0)=0$ and $h$ is not a $\cl^1$-diffeomorphism at
  $s=0$.  So we compute the limit:
  \begin{gather}
    \lim_{s\ra 0} \frac{\unga(s)}{s} = \lim_{t\ra 0 }
    \frac{\unga(h(t))} {h(t) } = \lim_{t\ra 0} \frac{\beta(t)}{h(t)}.
  \end{gather}
  Since $\beta(0) =0 $ and $h(0)=0$ we may apply de L'H\^{o}pital
  rule:
  \begin{equation}
    \label{eq:dunga}
    \lim_{s\ra 0} \frac{\unga(s)}{s} = \lim_{t \ra 0}
    \frac{\dot{\beta}(t)}{| \dot{\beta}(t)|} = \lim_{t\ra 0} \frac{
      \rho' + i \theta' \rho}{ \sqrt{ (\rho')^2 + i (\theta'\rho)^2}}
    e^{i\theta/m} \xi_k = e^{i\theta_0/m} \xi_k.  
  \end{equation}
  (Recall that $\rho'(0)=1$ and $\theta' \rho \ra 0$.)  This shows
  that $\unga $ is $\cl^1$ up to $s=0$. The last assertion is
  immediate from \eqref{eq:dunga}.
\end{proof}


\begin{lemma}\label{geodetiche-angolose}
  Let $\alfa: [0,\eps] \ra X$ and $\beta: [0, \eps] \ra X$ be segments
  with $\alfa(0)=\beta(0)=0$.  If $ \ango \bigl(\undalf(0), \undb(0)
  \bigr) < \pi/m $ then
  \begin{equation}
    \label{eq:meno-di-1}
    \lim_{t\ra 0} \frac{d(\alfa(t), \beta(t)) }{2t} =
    \sin \frac{  \ango ( \dalf(0), \db(0)) }{2}
    <1.
  \end{equation}
  In particular $\alfa * \beta^0$ is not minimising on any subinterval
  of $[-\eps, \eps]$ which contains $0$ as an interior point.
\end{lemma}
\begin{proof}
  By interchanging $\alfa$ and $\beta$ if necessary, we can assume
  that $ \undalf(0) = e^{i\eta_1} $ and $ \unbeta(0) = e^{i\eta_2} $
  with $0 \leq \eta_i < 2\pi$ and $0\leq \eta_2 - \eta_1 < \pi/m$.
  According to the previous Proposition $\dalf(0)=(\undalf^m(0), 0 ,
  \dotsc, 0)$ and $\db(0)=(\undb^m(0), 0 , \dotsc, 0)$.  Hence we can
  choose $\theta_i\in \R$ such that
  \begin{gather*}
    \begin{aligned}
      \dalf(0)&=(e^{i\theta_1}, 0 , \dotsc, 0) \\
      \db(0)& =(e^{i\theta_2}, 0 ,\dotsc, 0)
    \end{aligned}
    \qquad
    \begin{gathered}
      0\leq \theta_1 <2\pi \\
      0\leq \theta_2 - \theta_1 =\ango(\dalf(0), \db(0)) < \pi.
    \end{gathered}
  \end{gather*}
  We start by showing that
  \begin{gather}\label{eq:limite-angoloso}
    \lim_{t \to 0} \frac{|\beta_1(t) - \alfa_1(t)|}{2t} = \sin \frac{
      \ango ( \dalf(0), \db(0)) }{2}.
  \end{gather}
  We can find continuous function $\rho_\alfa, \rho_\beta,
  \theta_\alfa, \theta_\beta $ such that
  \begin{gather*}
    \begin{gathered}
      \alfa_1(t) =
      \rho_\alfa (t) e^{i \theta_\alfa(t)} \\
      \beta_1(t) = \rho_\beta (t) e^{i \theta_\beta(t)}
    \end{gathered}
    \qquad
    \begin{gathered}
      \theta_\alfa(0)=\theta_1 \\
      \theta_\beta(0)=\theta_2.
    \end{gathered}
  \end{gather*}
  And we know from \eqref{eq:rho-primo} that $
  \rho_\alfa(0)=\rho_\beta(0) =0$, $ \rho_\alfa'(0)=\rho'_\beta(0)
  =1$.  Therefore
  \begin{gather*}
    \lim_{t \to 0} \frac{|\beta_1(t) - \alfa_1(t)|}{2t} = \lim_{t \to
      0} \frac{1}{2t} \sqrt{ \rho^2 _\alfa + \rho^2_\beta -
      2\rho_\alfa \rho_\beta \cos (\theta_\beta - \theta_\alfa) } =
    \\
    = \demi{} \lim_{t\ra 0} \sqrt{ \Bigl (
      \frac{\rho_\alfa}{t}\Bigr)^2 + \Bigl (\frac{ \rho_\beta}{t}
      \Bigr)^2 - 2\frac{\rho_\alfa}{t} \frac{ \rho_\beta}{t} \cos
      (\theta_\beta - \theta_\alfa) }
    = \\
    =\demi \sqrt{(\rho_\alfa'(0))^2 + (\rho_\beta'(0))^2 -2
      \rho_\alfa'(0) \rho_\beta'(0) \cos( \theta_2-\theta_1)}=
    \\
    =\demi \sqrt{2 -2\cos (\theta_2-\theta_1)} = \sqrt{\frac{1 -\cos(
        \theta_2 - \theta_1)}{2}} =\sin \frac{\theta_2 - \theta_1
    }{2}.
  \end{gather*}
  Thus \eqref{eq:limite-angoloso} is proved.

  Next set $\theta_0 = (\theta_1 + \theta_2) / 2$.  Since $\theta_2 -
  \theta_1 < \pi$ both $\dalf_1(0) $ and $ \db_1(0)$ lie in the sector
  $ S(\theta_0, \pi/2)$.  By continuity there is $\delta_1>0$ such
  that
  \begin{gather*}
    \alfa_1((0, \delta_1]) \cup \beta_1((0, \delta_1]) \subset
    S(\theta_0, \pi/2).
  \end{gather*}
  For $j=0,1,\ldots{}, m-1$ set
  \begin{gather*}
    S_j := S\biggl (\frac{\theta_0+ 2\pi j}{m} , \frac{\pi}{2m} \biggr
    ).
  \end{gather*}
  Then $ u^\men \bigl ( S(\theta_0, \pi/2) \bigr ) = \sqcup
  _{j=0}^{m-1}S_j$.  Note that $ \alfa_1= \pi \alfa =\pi \phi
  \phi\meno \alfa = u \phi\meno \alfa$ and similarly $\beta_1 = u
  \phi\meno \beta$. Then $ \phi\meno \alfa(t) \in u^\men S(\theta_0,
  \pi/2)$ for $t\in (0, \delta_1)$ and by connectedness the image of
  $\phi\meno \circ \alfa $ must lie inside some component
  $S_j$. Similarly the image of $\phi\meno \circ \beta $ is entirely
  contained in some component $S_k$. Since the sectors $S_j$ and $S_k$
  are convex $\undalf(0) \in S_j$ and $\undb(0)\in S_k$ as well.  But
  $\eta_2 - \eta_1 \in (0, \pi/m)$ so $\undalf(0)$ and $\undb(0)$ in
  the same component of $u^\men \bigl ( S(\theta_0, \pi/2) \bigr)$.
  Hence $S_k=S_j$.  The restriction
  \begin{equation*}
    u_j:= u\restr{  S_j    }: S_j
    \ra S(\theta_0,\pi/2).
  \end{equation*}
  is a biholomorphism and
  \begin{gather*}
    \phi\meno\alfa(t) = u_j^\men (\alfa_1(t)) \qquad{}
    \phi\meno\beta(t) = u_j^\men (\beta_1(t)) .
  \end{gather*}
  Fix $t\in (0,\delta_1)$.  Since $S_0 \subset \C$ is a convex set the
  formula
  \begin{gather*}
    \la_t(s) = u_j^\men \bigl ( (1-s) \alfa_1(t) + s \beta_1(t) \bigr
    )
  \end{gather*}
  defines a path $ \la_t : [0,1]\ra \Delta$ and $ \mu_t := \phi \circ
  \la_t : [0,1]\ra X$ is a smooth path from $\alfa(t)$ to
  $\beta(t)$. Hence
  \begin{gather*}
    d(\alfa(t) , \beta(t)) \leq \Lo(\mu_t) = \int_0^1 \bigg| \frac{d
      \mu_t}{ds}(s) \bigg |_g\, ds.
  \end{gather*}
  Differentiating (in $s$) the identity $u \bigl( \la_t(s) \bigr) =
  (1-s) \alfa_1(t) + s \beta_1(t)$ we get
  \begin{gather*}
    \frac{d}{ds} u\bigl ( \la_t (s) \bigr ) \equiv \beta_1(t) -
    \alfa_1(t).
  \end{gather*}
  On the other hand
  \begin{gather*}
    \frac{d}{ds} u\bigl ( \la_t (s) \bigr ) = \frac{d}{ds} \bigl (
    \la_t (s) \bigr )^m = m \la^{m-1}_t(s) \frac{d \la_t}{ds}(s)
    \\
    \frac{d \la_t}{ds}(s)= \frac{\beta_1(t) - \alfa_1(t)} { m
      \la_t^{m-1}(s)}.
  \end{gather*}
  Using first \eqref{eq:def-Rz} and \eqref{eq:phi-primo-con-Rz} and
  next \eqref{eq:c0-metrico} and \eqref{eq:def-c0} we have
  \begin{gather*}
    \frac{d \mu_t}{ds}(s)= (\beta_1(t) - \alfa_1(t)) \bigl (e_1 +
    \la_t(s) R(\la_t(s))\bigr)
    \\
    \begin{aligned}
    \bigg | \frac{d \mu_t}{ds}(s) \bigg |_g& \leq \bigg | \frac{d
      \mu_t}{ds}(s) \bigg | ( 1 + c_0 |\mu_t(s)|) \\
&\leq |\beta_1(t) -
    \alfa_1(t)| (1 + c_0 |\la_t(s)|) ( 1 + c_0 |\mu_t(s)|) .  
    \end{aligned}
  \end{gather*}
  By \eqref{eq:c0-zeta} $ |\mu_t(s) | \leq c_0 |\la_t(s)|$ so
  \begin{gather*}
    \bigg | \frac{d \mu_t}{ds}(s) \bigg |_g \leq |\beta_1(t) -
    \alfa_1(t)| \bigl (1 + (c_0 + c_0^2) |\la_t(s)| + c_0^3 |\la_t(s)|
    \bigr) .
  \end{gather*}
  Moreover we have
  \begin{gather*}
    \begin{aligned}
      |\la_t(s)|^m &= |u(\la_t(s))|  =| (1-s) \alfa_1(t) + s\beta_1(s)| \leq \\
      & \leq (1-s) |\alfa(t)| + s |\beta(t)| \leq \\
      &\leq (1-s) d(0,\alfa(t)) + s \, d(0,\beta(t)) = t\\
    \end{aligned}
    \\
    |\la_t(s)| \leq t^\um.
  \end{gather*}
  So there is a constant $C>0$ such that
  \begin{gather}
    \nonumber{}
    (c_0 + c_0^2) |\la_t(s)| + c_0^3 |\la_t(s)| \leq  C t^\um\\
    \nonumber{} \bigg | \frac{d \mu_t}{ds}(s) \bigg |_x \leq
    |\beta_1(t) - \alfa_1(t)| (1 + C t^{\nicefrac{1}{m}})\\
    \nonumber{} d(\alfa(t) , \beta(t)) \leq |\beta_1(t) - \alfa_1(t)|
    (1 + C t^{\nicefrac{1}{m}}).
  \end{gather}
  This yields the upper bound
  \begin{gather*}
    \limsup_{t\ra 0} \frac{d(\alfa(t), \beta(t)) }{2t} \leq \lim_{t\to
      0} \frac{|\beta_1(t) - \alfa_1(t)| (1+Ct^\um)}{2t} = \lim_{t\ra
      0} \frac{|\beta_1(t) - \alfa_1(t)|}{2t} .
  \end{gather*}
  As for the lower bound, using \eqref{eq:distanze-asintotiche} we
  have
  \begin{gather*}
    \liminf_{t\ra 0} \frac{d(\alfa(t), \beta(t)) }{2t} \geq
    \liminf_{t\ra 0} \frac{|\alfa(t)- \beta(t)| }{2t} \cdot
    \liminf_{t\ra 0} \frac{d(\alfa(t), \beta(t)) }{|\alfa(t)-\beta(t)|} \geq \\
    \geq \liminf_{t\ra 0} \frac{|\alfa(t)- \beta(t)| }{2t} \geq
    \lim_{t\ra 0} \frac{|\alfa_1(t)- \beta_1(t)| }{2t}.
  \end{gather*}
  Thus using \eqref{eq:limite-angoloso} we finally compute the limit
  \begin{gather*}
    \lim_{t\ra 0} \frac{d(\alfa(t), \beta(t)) }{2t} = \lim_{t\ra 0}
    \frac{|\beta_1(t) - \alfa_1(t)|}{2t} = \sin \frac{ \ango (
      \dalf(0), \db(0)) }{2}.
  \end{gather*}
  Since $\ango ( \dalf(0), \db(0)) < \pi$ this completes the proof.
\end{proof}

\section{Uniqueness of geodesics}
\label{sec-uni}

\begin{lemma}\label{geodetiche-lisce}
  Let $\ga_1: [0,L_1 ] \ra X $ be a segment between two points $x,y\in
  \Xr$. If $\ga_2: [0,L_2 ] \ra X $ is a unit speed geodesic distinct
  from $\ga_1$ with $\ga_2(0)=x$ and $ \ga_2(t_2)=y$ for some $t_2 \in
  (0,L_2)$, then $\ga_2$ is not minimising beyond $t_2$, that is
  $d(\ga_2(0), \ga_2(t_t + \eps ) < t_2 + \eps$ for any $\eps>0$.
\end{lemma}
\begin{proof}
  If $\ga_2$ were minimising on $[0, t_2 + \eps]$, then $t_2=L_1$, the
  concatenation $\ga_1 * \ga_2\restr{[t_2, t_2+\eps]}$ would be a
  shortest path from $x$ to $\ga_2(t_2+\eps)$ and therefore would be
  smooth near $t_2$.  This would force $\ga_1 = \ga_2$.
\end{proof}
In the following we will repeatedly make use of the following
celebrated idea of Klingenberg (see \cite[Lemma 1]{klingenberg-59} or
\cite[Lemma 2.1.11(iii)]{klingenberg-libro}).
\begin{lemma}[Klingenberg]
  \label{klingenberg}
  Let $x$ be a point in a Riemannian manifold $(M,g)$ and let $v_1,
  v_2 \in U_xM$ be two distinct unit vectors such that $\ga^{v_1}$ and
  $\ga^{v_2}$ be defined and minimising on $[0,T]$. Assume that
  $\ga_{v_1}(T)=\ga_{v_2}(T)$, that $\dga_{v_1}(T) + \dga_{v_2}(T)
  \neq 0$ and that $\ga_{v_i}(T) $ is not a conjugate point of $x$
  along $\ga_{v_i}$. Then there are vectors $v_1',v_2'\in U_xX$
  arbitrarily close to $v_1 $ and $v_2$ respectively and such that the
  geodesics $\ga_{v_i'}$ are minimising on $[0,T']$ for some $T'<T$
  and $\ga_{v_1'}(T')=\ga_{v_2'}(T')$.
\end{lemma}

\begin{lemma} \label{simple-closed} Let $\ga_1 , \ga_2 : [0,L] \ra X$
  be segments with the same endpoints.  If $0 \not\in \ga_2([0,L))$
  then $\ga_1* \ga_2\no$ is a simple closed curve.
\end{lemma}
\begin{proof}
  Assume by contradiction that there are $t_1, t_2 \in (0,L)$ such
  that $\ga_1(t_1) =\ga_2(t_2)$.  Since $x=\ga_2(0) $ and
  $y=\ga_2(t_2)$ are regular points Lemma \ref{geodetiche-lisce}
  implies that $\ga_2$ is not minimising on $[0,L]$, contrary to the
  hypotheses.
\end{proof}
Since $X$ is a topological disc, by the Jordan separation theorem the
interior of a simple closed curve contained in $X$ is well defined and
is again a topological disc.  Fix on $\Xr$ the orientation given by
the complex structure.  If $\alfa : [0, L]\ra \Xr$ is a piecewise
smooth simple closed path in $X$ we say that it is \enf{positively
  oriented} if its interior lies on its left
\cite[p. 268]{do-carmo-surfaces}.  If $x\in \Xr$ and $u, v\in T_xX$
are two linearly independent vectors we let $\ango(u,v)$ denote the
\enf{unoriented angle} as before, while $\angor (u, v) $ denotes the
\enf{oriented angle}, which is defined by $\angor(u,v) = \ango(u,v)$
if $\{u,v\}$ is a positive basis of $T_xX$ and by $\angor(u,v) = -
\ango(u,v)$ otherwise. Equivalently, if $v=e^{i\theta} u$ with $\theta
\in (-\pi,\pi)$ then $\angor(u,v) = \theta$.  If $\alfa : [0, L]\ra
\Xr$ is a positively oriented piecewise smooth simple closed path and
$t\in (0,L)$ is a vertex that is not a cusp, the \enf{external angle}
at $\alfa(t)$ is defined as $\thext(t)=\angor(\dalf(t-), \dalf(t+))$,
and the \enf{interior angle} as $\thint(t)=\pi - \thext$. Note that
$\thext(t) \in (-\pi, \pi)$, while $\thint(t) \in (0, 2\pi)$
\cite[p.266ff]{do-carmo-surfaces}.

\begin{lemma}\label{GB}
  There is $r_1 > 0$ such that
  for any pair of segments $\ga_1 , \ga_2 : [0,L] \ra \Bd^*(0,r_1)$
  with the same endpoints
  $ \dga_1(0) \neq - \dga_2(0)$ and
  $ \dga_1(L) \neq - \dga_2(L)$.
  Moreover, if $
  \ga_1* \ga_2\no$ is positively oriented and $0$ does not lie in its
  interior, then
  the interior angles of $\ga_1* \ga_2\no$ at the two vertices are
  both smaller than $\pi$ and
  \begin{gather*}
    \angor\bigl(\dga_1(L) , \dga_2(L)\bigr ) < 0.
  \end{gather*}
\end{lemma}
\begin{proof}
%
%
  Using Lemma \ref{angoletto-3} we can find a $\delta>0$ with the
  following property: for any segment $\ga: [0,L]\ra \Bd(0,\delta)$
  with $\ga((0,L)) \subset \Xr$ we have
  \begin{gather}
    \label{eq:angoletto-55}
    \ango\bigl (\dga(s),\dga (s')\bigr ) < \deltino
  \end{gather}
  for any $s,s'\in [0,L]$, the angle being computed with respect to
  the Hermitian product $\langle \ , \, \rangle$.
  Let $\kl\in \R$ be defined as in \eqref{eq:def-K0}.
  By Wirtinger theorem \cite[p.159]{chirka} $\om$ is the volume form
  of $g\restr{\Xr}$.  By Lelong theorem \cite[p.173]{chirka} analytic
  sets have locally finite mass.  Hence there is an $r_1 \in (0,
  \delta) $ such that
  \begin{equation*}
    \vol\bigl (\Bd(0,\cuno^2 r_1^\um)\bigr ) =  
    \int_{\Bd(0,\cuno^2 r_1^\um)} \om < \frac{\pi}{1 + |\Kl| } 
  \end{equation*}
  Here $\cuno$ is the constant in Prop. \ref{param-hoelder}.  Let
  $\ga_1 , \ga_2 : [0,L] \ra \Bd^*(0,r_1)$ be a pair of
  segments 
  with the same endpoints.  Assume by contradiction that $ \dga_1(L) =
  - \dga_2(L)$. 
  Set $\alfa = \ga_1 * \ga_2^0$ and $w= \dga_1(L) = - \dga_2(L)$.  By
  \eqref{eq:angoletto-55}
  \begin{gather*}
    \ango \bigl ( \dalf(s) , w \bigr ) < \deltino \ \text{ i.e. }\
    \langle \dalf(s) , w\rangle >0
  \end{gather*}
  for any $s\in [0,2L]$. Hence
  \begin{gather*}
    \langle \alfa(2L) , w \rangle - \langle \alfa(0) , w \rangle =
    \int_0^{2L} \langle \dalf(s) , w \rangle ds >0.
  \end{gather*}
  In particular we would get $\ga_2(0) = \alfa(2L) \neq
  \alfa(0)=\ga_1(0)$ contrary to the hypothesis that the endpoints
  coincide.  This proves that $\dga_1(L) \neq - \dga_2(L)$. The same
  argument of course yields $\dga_1(0) \neq - \dga_2(0)$ as well.

  Next denote by $V$ be the interior of $\alpha$ and assume that
  $0\not \in V$ and that $\alfa$ is positively oriented.  By
  \eqref{eq:palletta-dischetto}
  \begin{equation*}
    \Bd(0,r_1) \subset U:=\phi(B(0, \cuno r_1^\um)) 
    \subset \Bd(0, \cuno^2 r_1^\um ) .
  \end{equation*}
  Since
%
  $U$ is a topological disc and $\partial V \subset U$, also $V\subset
  U \subset \Bd(0, \cuno^2 r_1^\um ) $.
  Since $\overline{V}\subset \Xr$ Gauss--Bonnet theorem applies and we
  get
  \begin{equation*}
    \thint(0) + \thint(L) = \int _{V} K\om \leq \Kl\cdot
    \vol\bigl(\Bd(0,\cuno^2r^\um_1)\bigr) < \pi.
  \end{equation*}
  Thus $\thint(0), \thint(L) \in [0,\pi)$.  To prove the last
  assertion set $\theta = \angor\bigl ( \dga_1(L), \dga_2(L)\bigr )$.
  Since $\ga_1$ and $\ga_2$ are distinct geodesics $\theta \neq 0$.
  It is easy to check that
  \begin{gather*}
    \thint(L) =
    \begin{cases}
      2\pi -\theta  & \text {if } \theta \in (0, \pi) \\
      - \theta & \text {if } \theta \in (- \pi, 0).
    \end{cases}
  \end{gather*}
  Since $\thint(L) \in (0, \pi)$, $\theta \in (-\pi, 0)$.
%
%
%
%
%
%
%
%
\end{proof}
Let $\delta >0$ be such that $\Bd(0,\delta) \compsubset{} X$.  Put
\begin{equation*}
  r_2=\demi \min\biggl \{ \frac{\pi}{\sqrt{\Kl}},
  {\delta}, r_1 
  \biggr \} 
\end{equation*}
where $\pi/\sqrt{\Kl} = + \infty $ if $\Kl\leq 0$.

\begin{prop}\label{surmolotto}
  For any $x\in \Bd(0,r_2)$, $B_x\bigl (0, d(x,0)\bigr ) \subset
  \Dom_x$, $\exp_x$ has no critical points on $\bx(0,r_2)\cap \Dom_x$
  and $\exp_x \bx(0,$ $d(x,0)) = \Bd\bigl(x,d(0,x)\bigr)$.
\end{prop}
\begin{proof}
  Let $x\in \Bd(0,r_2)$ and $v\in U_xX $.
  Set $r=d(x,0)$ and let $T_v$ be as in \eqref{eq:def-Tv}.  Assume by
  contradiction that $T_v<r$ and set $\eps= (r-T_v)/2>0$.  For any
  $t\in [0,T_v)$
  \begin{gather*}
    d(\ga^v(t), 0) \leq d(\ga^v(t), x) + d(x,0) \leq t + r < 2r \leq
    2r_2
    \leq \delta \\
    d(\ga^v(t), 0) \geq | d(\ga(t),x) - d(x, 0)| \geq r-t \geq r-T_v >
    \eps.
  \end{gather*}
  So $ \ga^v([0,T_v)) $ is contained in $ Q:=\overline { \Bd(0,
    \delta)}\setminus \Bd(0, \eps )$ and $t\mapsto \dga^v(t)$ is a
  trajectory of the geodesic flow contained in the compact set
  $\{(y,w)\in TX_\reg: y\in Q, |w| =1\}$. This contradicts the
  maximality of $T_v$. Therefore $T_v\geq r$ and $\bx(0,r) \subset
  \Dom_x$.  Since $K \leq \Kl$ on $\Xr$ and $r_2\leq \pi/\sqrt{\Kl}$
  Rauch theorem \cite[p.215]{do-carmo} implies that for any $v\in
  U_xX$ the geodesic $\ga^v$ has no conjugate points on
  $[0,\min\{T_v,r_2\})$.  Therefore $\exp_x$ is a local diffeomorphism
  on $\bx(0,r_2)\cap \Dom_x$ \cite[p.114]{do-carmo}. This proves the
  second claim.  The inclusion $\exp_x \bx(0, r) \subset \Bd(x,r)$ is
  obvious.  On the other hand if $y \in \Bd(x,r)$ let $\ga : [0,
  d(x,y)] \ra X$ be a segment from $x$ to $y$.  By the triangle
  inequality
  $\ga$ is contained in $\Xr$ so $\ga(t)= \exp_xtv$ for some $v\in
  U_xX$.  Then $y=\exp_x d(x,y)v \in \exp_x \Bd(0,r)$. This proves
  that $\exp_x\bx(0,r) = \Bd(x,r)$.
\end{proof}
For $x\in \Bd(0,r_2)$ define $c_x: U_xX \ra (0,r_2]$ by
\begin{equation}
  \label{eq:def-di-c} 
  c_x(v) = \sup \{ t \in (0,\min\{T_v, r_2\}): 
  \ga^v 
  \text { is minimising on } [0,t]\}
\end{equation}
and put $\cc(x) = \inf _{\ux} c_x$. If $\ga^v$ is a segment from $x$
to $0$ then $c_x(v)=d(x,0)$, so $\cc_x \leq d(x,0)$.  In the next two
lemmata we adapt to our situation arguments that are classical in the
study of the cut locus of a complete Riemannian manifold, see
e.g. \cite[p.102]{sakai-Riemannian}.  For the reader's convenience we
provide all the details.
\begin{lemma}
  \label{caratterizzo-c}
  Let $x\in \Bd(0,r_2)$, $v\in\ux$ and $T\in (0, \min\{T_v, r_2\})$.
  Then $T=c_x(v)$ iff $\ga^v$ is minimising on $[0,T]$ and there is
  another segment $\ga \neq \ga^v$ between $x$ and $\ga^v(T)$.  If
  $d(x,0) + d(0, \ga^v(T)) > T$ then $\ga$ lies entirely in $\Xr$, so
  $\ga^u=\ga$ for some $u\in \ux$, $u \neq v$.  In particular this
  happens if $T< d(x,0)$.
\end{lemma}
\begin{proof} Put $y=\ga^v(T)\in \Xr$ and assume $T=c_x(v)$. Then
  $\ga^v$ is minimising on $[0,t]$ for any $t<T$, so also on $[0,T]$.
  Since it is not minimising after $T$, we may choose a sequence $t_n
  \searrow T$ such that $\ga^v$ is never minimising on $[0, t_n]$.
  Put $ y_n= \ga^v(t_n)$ and $ s_n=d(x, y_n)$. Then $s_n < t_n$ and
  $s_n \ra T$. Let $\ga_n : [0, s_n ] \ra X $ be a segment from $x$ to
  $y_n$.  By Ascoli-Arzel\`a Theorem and Cor. \ref{semicont} we can
  extract a subsequence converging to a segment $\ga: [0,L]\ra X$ from
  $x$ to $y$.  If $\ga=\ga^v$, then $\ga_n$ is contained in $\Xr$ for
  large $n$, so $\ga_n=\ga^{v_n}$ for some $v_n \in U_xX$ and $v_n \ra
  v$. But then any neighbourhood of $Tv$ in $T_xX$ contains a pair of
  distinct points $s_n v_n \neq t_n v$ that are mapped by $\exp_x$ to
  the same point $y_n \in \Xr$. Since $Tv \in \bx(0,r_2)\cap \Dom_x$
  this contradicts Prop. \ref{surmolotto}.  Therefore $\ga\neq \ga^v$.
  This proves necessity of the condition.  Sufficiency follows
  directly from
  Lemma \ref{geodetiche-lisce}.  The remaining assertions are trivial.
\end{proof}

\begin{lemma} \label{c-semic} For $x\in \Bd(0,r_2)$ the function $c_x
  $ is lower semicontinuous. In particular the minimum $\cc_x$ is
  attained.
\end{lemma}
\begin{proof}
  Let $v_n\in\ux$ be a sequence such that $v_n\to v$.
  Set $T:= \liminf_{n\to \infty} c_x(v_n)$.  We wish to prove that $
  c_x(v) \leq T$.  If $T=r_2$ this is obvious from the definition
  \eqref{eq:def-di-c}.
  Assume instead that $T < r_2$.  Passing to a subsequence we can
  assume that $T_n:=c_x(v_n) < r_2$ and $T_n \to T$.  By the theorem
  of Ascoli-Arzel\`a the segments $\ga^{v_n} \restr{[0,T_n]}$ converge
  to
  a segment $\alfa : [0,T] \ra X$ 
  and $\alfa(t)=\ga^v(t)$ for $t\in[0,T_v)$. If there is $\tau\in
  (0,T]$ such that $\alfa(\tau)=0$ then $c_x(v) \leq T_v \leq \tau
  \leq T$ and we are done.  Otherwise $\alfa ([0,T]) \subset \Xr$, so
  $\ga^v=\alfa$ is minimising on $[0,T]$ and $T < T_v$.
%
%
%
  For $n$ large $\ga^{v_n}([0,T_n] \subset \Xr$ as well.  Hence
  $c_x(v_n) < \min \{T_v, r_2\}$.  By the previous lemma there are
  segments $\ga_n\neq \ga^{v_n}$ from $x$ to $\ga_{v_n}(T_n)$. Again
  by the theorem of Ascoli-Arzel\`a we can assume, by passing to a
  subsequence, that $\ga_n $ converge to a segment $ \beta$ from $x$
  to $\ga^v(T)$.  If $\beta $ passes through $0$ then $\beta \neq
  \ga^v$ and 
  $T=c_x(v)$ by the previous lemma.
  If $\beta$ is contained in $\Xr$, the same is true of $\ga_n$ for
  large $n$. Write $\ga_n=\ga^{u_n}$ and extract a subsequence so that
  $u_n \to u$. Clearly $\beta=\ga^u$. 
%
%
  If $u=v$, any neighbourhood of $Tv$ would contain two distinct
  vectors $T_nv_n \neq T_n u_n$ with the same image through
  $\exp_x$. Since $T < r_2$ this possibility is ruled out by
  Prop. \ref{surmolotto}. Therefore $u\neq v$ and the previous lemma
  implies that $c_x(v)=c_x(u)=T$. 
\end{proof}

\begin{prop} \label{inietto} If $x\in \Bd(0,r_2)$ then
  $\cc_x=\inj{x}=d(x,0)$.
\end{prop}
\begin{proof}
  Let $x\in \Bd(0, r_2)$ and $r=d(x,0)$.
  First of all we prove that $\exp_x$ is injective on $\bx(0,\cc_x)$.
  In fact let $w_1, w_2 \in \bx(0,\cc_x) $ be such that
  $\exp_x(w_1)=\exp_x(w_2)$. Write $w_i=t_iv_i$ with $|v_i|=1$.  Since
  $ t_i =|w_i|<\cc_x \leq c(v_i)$ the geodesics $\ga^{v_i}$ are
  minimising on $[0,t_i]$. Therefore $t_1=d(x,\exp_x(w_i))=t_2$.  If
  $v_1\neq v_2$, Lemma \ref{caratterizzo-c} would imply that
  $t_1=c(v_1)$, but this is impossible since $t_1 < \cc_x$.  So
  $v_1=v_2 $ and $w_1=w_2$. This proves that $\exp_x$ is injective on
  $\Bd(0,\cc_x)$.  Next we prove that $\cc_x=r$.  We already know that
  $\cc_x \leq r$.
  Assume by contradiction that $T:=\cc_x < r$.  By Lemma \ref{c-semic}
  there are $u\neq v\in \ux$ such that $\ga^u(T) = \ga^v(T)$. Since
  $T<r$ Prop. \ref{surmolotto} ensures that $\exp_x$ is a
  diffeomorphism on appropriate neighbourhoods of $Tu$ and $Tv$ in
  $T_xX$.  By Lemma \ref{klingenberg} we conclude that $\dga^u(T) = -
  \dga^v(T)$.  But this is impossible by Lemma \ref{GB}.
  Therefore $\cc_x=r$ and $\exp_x$ is injective on $\bx(0,r)$. Hence
  $\exp_x$ is a diffeomorphism of $B_x(0,r)$ onto $\Bd(x,r)$. In
  particular $\inj{x} \geq \cc_x \geq r$.  The reverse inequality is
  proven in Lemma \ref{inj-distanza-da-sing}.
\end{proof}

\begin{prop}\label{unicita-regolare}
  There is $r_3\in (0,r_2/2)$ such that if $\alfa, \beta : [0,T]\ra
  \Xr$ are distinct segments with the same endpoints $x, y \in
  \Bd^*(0,r_3) $, then $0$ lies in the interior of $\alfa*\beta^0$.
\end{prop}
\begin{proof}
  Set
  \begin{equation*}
    r_3=
    \Bigl (\frac{r_2}{2\cuno^{2}} \Bigr )^m 
  \end{equation*}
  where $\cuno$ is the constant in \eqref{eq:hold-phi}.  By
  \eqref{eq:palletta-dischetto}
  \begin{equation}
    \Bd(0,r_3) \subset U:=\phi(B(0, \cuno r_3^\um)) 
    \subset \Bd(0, \cuno^2 r_3^\um ) \subset\Bd(0,r_2/2)    
  \end{equation}
  and $U$ is a topological disc.  
  Let $\alfa$ and $\beta $ be as above and set $x=\alfa(0)=\beta(0),
  y=\alfa(T)=\beta(T)$.  Since $x, y \in \Bd(0,r_3) \subset
  \Bd(0,r_2/2)$, $\alfa$ and $\beta$ lie in $\Bd(0, r_2) \subset
  \Bd(0,r_1)$.  By Lemma \ref{simple-closed} $\alfa*\beta^0$ is a
  simple closed curve.  Denote by $V$ its interior and
  assume by contradiction that $0\not \in V$.  Since $ \partial V
  \subset U$ also $\overline{V}\subset U\subset \Bd(0,r_2/2)$ and
  $\operatorname{diam} \overline{V} < r_2$.
  In particular $T < r_2$.  By interchanging $\alfa$ and $\beta$ we
  can assume that $V$ lies on the left of $\alfa*\beta^0$.  Set
  $u_0=\dalf(0) \in U_xX$ and chose $\theta_0\in (0,2\pi)$ so that
  $\dot{\beta}(0)=e^{i\theta_0}u_0 $.  By hypothesis $\theta_0>0$.
  Denote by $E$ the set of unit vectors $v\in \ux$ of the form
  $v=e^{i\theta} u_0 $ with $\theta\in [0,\theta_0]$ and by $\intE$
  the subset of those with $\theta \in (0,\theta_0)$.  By Lemma
  \ref{c-semic} the function $c_x$ has a minimum on $E$.
  \begin{figure}[h]
    \begin{center}
      \includegraphics{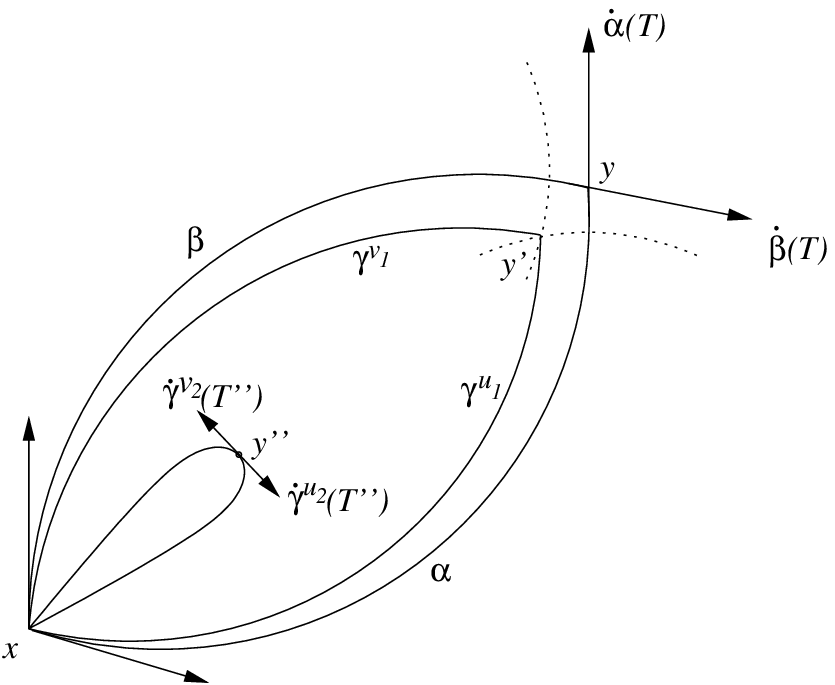}
    \end {center}
    \caption{}
    \label{baffi}
  \end{figure}
%
%
  We claim that the minimum point lies in $\inte E$.  By the
  hypotheses $c_x\bigl (\dalf(0) \bigr) = c_x \bigl(\db(0)\bigr) = T <
  r_2$.  Lemma \ref{GB} implies that $\theta_0< \pi$ and $\angor\bigl
  ( \dalf(T),\dot{\beta}(T) \bigr ) \in (- \pi,0)$ (see
  Fig. \ref{baffi}).  By Klingenberg Lemma \ref{klingenberg} there are
  vectors $u_1$ and $v_1$ arbitrarily close to $\dalf(0)$ and
  $\dot{\beta}(0)$ respectively, such that $T'=c_x(u_1)=c_x(v_1) < T$
  and $\ga^{u_1}(T')=\ga^{v_1}(T')$.  Since $\angor\bigl ( \dalf(T),
  \db(T) \bigr ) \in (-\pi, 0)$ the point
  $y'=\ga^{u_1}(T')=\ga^{v_1}(T')$ belongs to $V$.  Therefore
  $\ga^{u_1}(t)$ and $\ga^{v_1}(t)$ lie inside $V$ for any $t\in (0,
  T']$: otherwise they would meet either $\alfa$ or $\beta$ at an
  interior point, which is forbidden by Lemma
  \ref{geodetiche-lisce}. This shows that $u_1, v_1 \in \inte E$ and
  that $\dalf(0)$ and $\db(0)$ are not local minima of $c_x\restr{E}$
  and that the minimum of $c_x$ on $E$ must be attained at some point
  $u_2 \in \inte E$. Set $T''=c_x(u_2)=\min_E c_x$ and
  $y''=\ga^{u_2}(T'')$. Since $u_2 \in \inte E$ and $T''< T$, the
  point $y''$ belongs to $ V$, so $\ga^{u_2}(t) \in V$ for any $t\in
  (0,T'']$ (use again Lemma \ref{geodetiche-lisce}).  By Lemma
  \ref{caratterizzo-c} there is a segment $\ga \neq \ga^{u_2}$ between
  $x$ and $y''$ and, again by Lemma \ref{geodetiche-lisce}, it is
  contained in $V$ as well.  So $\ga=\ga^{v_2}$ for some $v_2\in \inte
  E$ and $c_x(v_2)=d(x, y'')= c_x(u_2)$.
%
%
%
  Since $y''\in V \subset \Bd(0, r_2/2)$ both $\ga^{u_2}$ and
  $\ga^{v_2}$ are contained in $\Bd^*(0, r_2 ) \subset \Bd^*(0, r_1)
  $. Hence by Lemma \ref{GB} $ \ga^{u_2}(T'') \neq -
  \ga^{v_2}(T'')$. But then we can apply again Klingenberg lemma to
  get a pair of nearby vectors with $c_x$ strictly smaller than
  $T''$. Since $T''$ is the minimum this yields the desired
  contradiction.
\end{proof}

\begin{teo}
  \label{unicita}
  There is $\rqua \in (0, r_3)$ such that for any $x\in \Bd(0,\rqua)$
  there is a unique segment from $x$ to $0$.
\end{teo}
\begin{proof}
  Set
  \begin{equation*}
    \rqua=  \Bigl (\frac{r_3}{\cuno^{2}} \Bigr )^m 
  \end{equation*}
  where $\cuno$ is the constant in \eqref{eq:hold-phi}.  By
  \eqref{eq:palletta-dischetto}
  \begin{equation}
    \label{eq:palla-disco-palla}
    \Bd(0,\rqua) \subset U:=\phi(B(0, \cuno \rqua^\um)) 
    \subset \Bd(0, \cuno^2 \rqua^\um ) \subset\Bd(0,r_3)    
  \end{equation}
  and $U$ is a topological disc.  Fix $x\in \Bd(0,\rqua)$ and assume
  by contradiction that there are two distinct segments $\alfa, \beta
  : [0, r] \ra X $ from $x$ and $0$.  By Lemma \ref{simple-closed}
  $\alfa*\beta^0$ is a simple closed curve.  Let $V$ be the interior
  of $\alfa*\beta\no$.  Since $\partial V \subset U$ also
  $\overline{V} \subset U \subset \Bd(0,r_3) \subset \Bd(0,r_2/2)$.
%
%
  Assume that $V$ lies on the left of $\alfa*\beta^0$ and set
  $u_0=\dalf(0) \in U_xX$, $\dot{\beta}(0)=e^{i\theta_0}u_0 $ with
  $\theta_0 \in(0,2\pi) $.  Denote by $E$ be the set of $v\in \ux$ of
  the form $v=e^{i\theta} u_0 $ with $\theta\in [0,\theta_0]$ and by
  $\intE$ the subset of those with $\theta \in (0,\theta_0)$.  If
  $v\in\intE$ then $\ga^v(t) \in V$ for small positive $t$.  Set
  $r=d(x,0)$.  By
  Prop. \ref{inietto} $\ga^v$ is defined and minimising on $[0,r)$.
  Let $[0,T_v)$ be the maximal interval of definition of $\ga^v$. If
  $\ga^v((0,T_v))$ were not contained in $V$, there would be a minimal
  time $t_0$ such that $\ga^v(t_0) \in \alfa((0, L]) \cup
  \beta((0,L])$.  By Lemma \ref{geodetiche-lisce} this would imply
  that $t_0> c_x(v)$, so there would be a point $y=\ga^v(c_x(v)) \in V
  $ that is reached by two distinct segments starting from $x$. But
  this is impossible because of Prop. \ref{unicita-regolare} because
  $x, y \in \Bd(0,r_3)$. Therefore $\ga^v((0,T_v))$ has to be
  contained in $V$.  This implies that $c_x(v) \leq \diam V < r_2$.
  If $c_x(v) < T_v$ Lemma \ref{caratterizzo-c} would give again a pair
  of distinct segments with the same endpoints $x, \ga^v(c_x(v)) \in V
  \subset \Bd(0,r_3)$, thus contradicting
  Prop. \ref{unicita-regolare}. So $c_x(v)=T_v$.  Now $\ga^v$ is
  minimising hence Lipschitz on $[0,c_x(v))$ and therefore extends
  continuously to $[0,c_x(v)]$. The only possibility is that
  $\ga^v(c_x(v))=0$ and $c_x(v)=d(x,0) =r$.
  Let $S$ be the set of vectors $v\in T_xX$ of the form $v=\rho
  e^{i\theta }v_1$ with $\rho \in (0, r)$ and $\theta \in (0,
  \theta_0)$.  We have just proved that the map
  \begin{gather*}
    F: \overline{S} \ra \overline{V} \qquad F(w) =
    \begin{cases}
      \exp_x(w) & \text { if } |w|< r\\
      0 & \text { if } |w|=r
    \end{cases}
  \end{gather*}
  is continuous. Both $\overline{S}$ and $\overline{V}$ are
  topological discs, $F(\partial S ) \subset \partial V$ and
  $F\restr{\partial S } : \partial S \ra
  \partial V$ has degree 1 so it is not homotopic to a
  constant. Therefore $F$ must be onto, $\exp_x(S)=V$ and $V \subset
  \Bd(x,r)$.  Now we look at our configuration of geodesics from the
  point of view of $\Delta$ as in \S \ref{regu-delta}.  Set $\ga_1 =
  \alfa^0$ and $\ga_2=\beta^0$ and let $\unga_i: [0, L_i]\ra \Delta$
  be as in Def. \ref{def-unga}.  By Prop.  \ref{rego-in-D} these
  $\unga_i$ are $\cl^1$ paths on $[0,L_i]$.  For small $s$, each of
  them intersects the circle $Z_s=\{\z \in \Delta: |\z|=s\}$ at
  exactly one point $p_i(s)$.  Let $t_i(s) \in (0, r)$ be such that
  $p_i(s) = \phi \meno (\ga_i(t_i(s)))$.  The functions $t_i :[0,\eps
  ) \ra [0,r]$ are continuous, strictly decreasing in a neighbourhood
  of $0$ and such that $t_i(0)=0$. Since $\ga_1$ and $\ga_2$ do not
  intersect except at their endpoints $p_1(s) \neq p_2(s)$. Therefore
  the circle $Z_s$ is cut by $p_1(s)$ and $p_2(s)$ in exactly two
  arcs. One of them lies in $\phi\meno(V)$ the other outside of
  it. Denote by $\beta_s : [0,1] \ra \Delta$ a $\cl^1$ parametrisation
  of the former.  Then $\alfa_s:= \phi \circ \beta_s$ is a path of
  length $ L(\alfa_s) \leq ||d\phi||_{\infty} L(\beta_s) \leq 2\pi c_0
  s =C s$ lying in $\exp_x(\bx(0,r))$ and connecting $\ga_1(t_1(s))=
  \exp_x((r-t_1(s)) u_0)$ to $\ga_2(t_2(s))=
  \exp_x((r-t_2(s))e^{i\theta_0} u_0)$.  On the other hand 
  \begin{gather*}
    r < r_2 \leq \frac{\pi}{2\sqrt{\Kl}}
  \end{gather*}
  $K\leq \Kl$ on $\Bd(x,r)$ and $\exp_x$ is a diffeomorphism of
  $\bx(0,r) \subset T_xX $ onto $\Bd(x,r)$. Therefore a classical
  corollary to Rauch theorem \cite[Prop. 2.5 p.218]{do-carmo} ensures
  that $ L(\alfa_s) $ is bounded from below by some positive constant
  depending only on $\theta_0$ and $\Kl$.  This yields the
  contradiction and shows that the segments $\alfa$ and $\beta$
  coincide.
\end{proof}

\begin{lemma}\label{lemma-r5}
  There is $\rcin \in (0, \rqua/3)$ such for any segment $\ga: [0,L]
  \ra \Bd(0,3\rcin) $ with $\ga((0,L)) \subset \Xr$ and for any
  $s,s\in [0,L]$
  \begin{gather}
    \label{eq:angoletto-5}
    \ango(\pi(\dga(s)),\pi(\dga (s'))) < \deltone.
  \end{gather}
\end{lemma}
\begin{proof}
  By Lemma \ref{angoletto-4} there is $\delta>0$ such that
  \eqref{eq:angoletto-5} holds for any segment $\ga: [0,L]\ra
  \Bd^*(0,\delta)$.  Set $ \rcin= \min\{\delta, \rqua/2\}$.
\end{proof}

If $\ga : [a,b]\ra \C^*$ is a continuous path define its winding
number by
\begin{equation}
  \label{eq:def-W}
  W(\ga) = \Re \int_\ga \frac{dz}{z}.
\end{equation}
For a non-closed path $W(\ga) \in \R$.  The winding number $W(\ga)$
depends only on the homotopy class of $\ga$ with fixed endpoints. If
$\ga(t) = \rho(t)e^{2\pi i \theta(t)}$ with $\theta \in C^0([a,b])$
then
\begin{equation}
  \label{eq:W-theta}
  W(\ga)  = \theta(b) - \theta(a).
\end{equation}

\begin{lemma}\label{non-giro}
  If $\ga : [0,L]\ra \Bd^*(0,3\rcin)$ is a segment then $W(\pi\circ
  \ga) <1$.
\end{lemma}
\begin{proof}
  Set $\alfa=\pi\circ \ga$ and write $ \alfa(t) = \rho(t)
  e^{i2\pi\theta(t)} $ with $\theta \in C^0([0,L])$. Then
  $W(\alfa)=\theta(L) - \theta(0)$.  Assume by contradiction that
  $W(\alfa)\geq 1$.  Pick $t_0\in [0,L]$ such that
  $\theta(t_0)-\theta(0)=1$ and let $\chi :[0,t_0+1]\ra \Delta$ be
  defined by
  \begin{gather*}
    \chi(t) =
    \begin{cases}
      \alfa(t) & t\in [0,t_0]\\
      \alfa(t_0) + (t-t_0) ( \alfa(0)- \alfa(t_0)) & t\in [t_0, t_0+1]
      .
    \end{cases}
  \end{gather*}
  The second piece of $\chi$ is a parametrisation of the segment from
  $\alfa(t_0) $ to $\alfa(0)$. Since $\theta(t_0)-\theta(0)=1$, $\chi$
  is a loop that avoids the origin and has winding number $1$, so its
  homotopy class is a generator of $\pi_1(\Delta^*,\alfa(0))$.  Set
  \begin{gather*}
    v= \frac{\alfa(0) -\alfa(t_0)} {|\alfa(0) -\alfa(t_0)|} \qquad
    \begin{gathered}
      u_1(t)= \chi(t)\cdot v \\
      u_2(t)= \chi(t)\cdot Jv
    \end{gathered}
  \end{gather*}
  \begin{figure}[h]
    \begin{center}
      \includegraphics{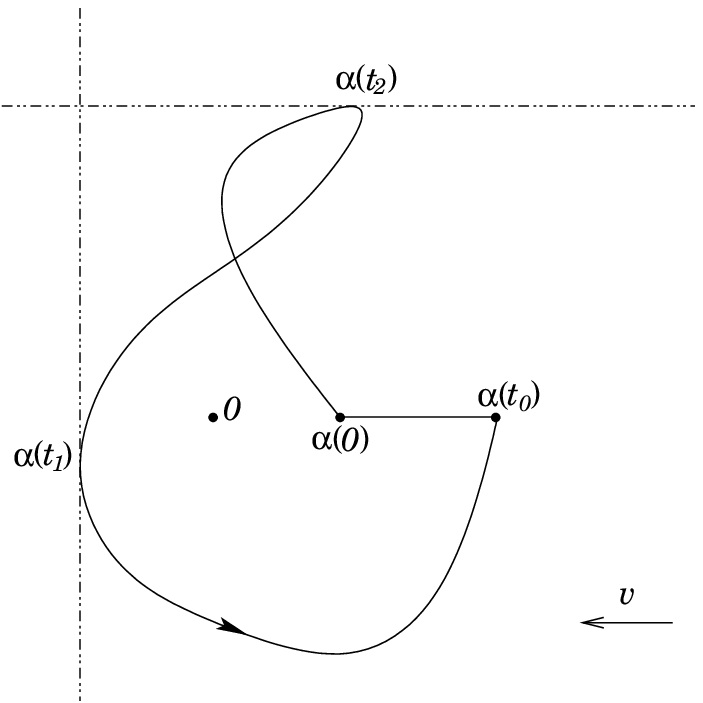}
    \end {center}
    \caption{}
    \label{tucano}
  \end{figure}
  ($J$ is the complex structure on $\C$.)  Since $W(\chi) =1$ both
  functions $u_1$ and $u_2$ have positive maximum, so their maximum
  points $t_1$ and $t_2$ belong to $(0, t_0)$.  Therefore $\dalf(t_2)
  = \pm v$, $\dalf(t_1) = \pm v$ and $\ango (\dalf(t_1), \dalf(t_2))
  \geq \pi/2$ (see Fig. \ref{tucano}). But $\alfa$ is the projection
  of the segment $\ga$. This contradicts \eqref{eq:angoletto-5} and
  proves the lemma.
\end{proof}

\begin{teo}
  \label{unicita-giro}
  For any $x, y \in \Bd^*(0,\rcin)$ there is at most one segment from
  $x$ to $y$ avoiding $0$.
\end{teo}
\begin{proof}
  Let $\ga_1, \ga_2 : [0,L] \ra X$ be two segments from $x$ to
  $y$. Both $\ga_1$ and $\ga_2$ are contained in $\Bd(0,3\rcin)$.
  Assume by contradiction that the two segments are distinct and both
  lie in $\Xr$.  By Lemma \ref{simple-closed} $\ga=\ga_1*\ga_2^0$ is a
  Jordan curve
  and by Prop. \ref{unicita-regolare} the origin lies in the interior
  of $\ga$.  Hence $[\ga]$ is a generator of $\pi_1(\Xr, \ga(0))$.
  Set $\alfa_i=\pi\circ \ga_i$, $\alfa=\pi\circ \ga=\alfa_1 *
  \alfa_2^0$.  Since $\pi: \Xr\ra \Delta^*$ is a degree $m$ unramified
  covering, the loop $\alfa$ has winding number $m$. Therefore either
  $\alfa_1$ or $\alfa^0_2$ has winding number at least 1. Nevertheless
  this is impossible by Lemma \ref{non-giro}.
\end{proof}

\section{Convexity}
\label{convex-section}

For $x\in \Bd^*(0, \rcin)$ denote by
\begin{equation}
  \nonumber{}
  \ga_x:
  [0, d(0,x)] \ra X    
\end{equation}
the unique segment from $0$ to $x$.  Define three maps
\begin{gather}
  \begin{aligned}
    &     \cobra: \Bd^*(0, \rcin) \ra S^1\times \{0\}\subset  \C^n \\
    &     \uncobra: \Bd^*(0, \rcin) \ra S^1\\
    & \camel: \Bd^*(0,\rcin) \times [0,1] \ra X
  \end{aligned}
  \qquad
  \begin{aligned}
    &     \cobra(x)=\dga_x(0)\\
    &        \uncobra(x)=\dunga_x{}(0) \\
    & \camel (x,t) = \ga_x \bigl ( td(0,x) \bigr ).
  \end{aligned}
  \label{eq:def-camel}
\end{gather}
%
%
$\cobra$ takes values in $ S^1\times \{0\}$ because $C_0X = \C \times
\{0\}$ and $g_x = \langle \ , \, \rangle$.
\begin{prop}\label{F-cont}
  The maps $\cobra$, $\uncobra$ and $\camel$ are continuous and
  $\cobra = (u \circ \uncobra, 0,\ldots, 0)$.
\end{prop}
\begin{proof}
  Assume $x_n \to x$ and set $\ga^n=\ga_{x_n}$.  By Theorem
  \ref{W-Hoelder} $||\dot{\ga}^n||_{\cl^{0,\um}} \leq \ctre$. By the
  Ascoli-Arzel\`a theorem there is a subsequence, still denoted by
  $\ga^n$, that converges in the $\cl^1$-topology to the unique
  segment $\ga_x $ from $0$ to $x$. In particular $\cobra(x_n^*) =
  \dga^n(0) \to \dga(0)= \cobra(x)$.  This shows that $\cobra$ is
  continuous.  That $\cobra = u \circ \uncobra$ was already proved in
  Prop.  \ref{rego-in-D}.  If $ \dga(0) =( e^{i\theta_0}, \ldots{},
  0)$, pick $t_0 \in (0, d(x,0))$ sufficiently close to $0$ that
  $\ga_1(t_0) \in S(\theta_0, \pi/2)$.  Denote by $S_1, \ldots, S_m$
  the connected components of $u\meno (S(\theta_0, \pi/2))$ and assume
  that $\phi\meno \ga(t_0)\in S_j$.  Then $\phi\meno \ga((0, t_0]) $
  is entirely contained in $ S_j$. Since $S_j$ is convex it follows
  that $\dunga(0) \in S_j$.  As $\ga^n \to \ga$ uniformly and
  $\phi\meno$ is H\"older (Prop. \ref{param-hoelder}) also $\phi\meno
  \ga^n(t_0)\in S_j$ and $\dunga^n(0) \in S_j$ for large $n$.  The map
  $u_j=u\restr{S_j}: S_j \ra S$ is a homeomorphism and
  $u(\dunga^n(0) )= \dga_1^n (0)$ and $u(\dunga(0) )= \dga_1 (0)$.
  Therefore $\dunga^n(0) = u_j \meno \dga_1^n (0) \to u_j \meno \dga_1
  (0) = \dunga(0) $.
%
  This proves that $\uncobra$ is continuous.  Finally, if $x_n \to x$
  and $t_n \to t$, by passing to a subsequence we can assume that
  $\ga^n =\ga_{x_n}\to \ga_x$ uniformly.
  Then clearly $\camel (x_n,t_n)=\ga^n\bigl (t_n d(x_n, 0)\bigr) \to
  \ga_x\bigl (t d(x,0) \bigr)= \camel (x,t)$.  This proves that the
  third map $\camel $ is continuous.
\end{proof}

\begin{prop}\label{piccolo-Whitehead}
  For any pair of points $x,y \in \Bd^*(0,\rcin)$ with
  \begin{gather}
    \label{eq:angolo-piccolo}
    \ango (\uncobra(x), \uncobra(y)) < \frac{\pi}{m}
  \end{gather}
  there is a unique segment $ \alfa_{x,y} : [0,d(x,y)] \ra X$ such
  that $\alfa_{x,y}(0)=x$ and $\alfa_{x,y}(d(x,y))=y$. This segment
  lies entirely in $\Xr$. If $\alfa_{x,y}(t) = \exp_x tv$, then $
  d(x,y) < c_x(v)$. Finally the map $(x,y,t) \mapsto \alfa_{x,y}(t)$
  is continuous.
\end{prop}
\begin{proof}
  Since $\ango(\uncobra(x), \uncobra(y)) < \pi/m$ it follows from
  Lemma \ref{geodetiche-angolose} that the path $ \ga_x * \ga_y^0 $ is
  not minimising.  If $\ga_1$ and $\ga_2$ were two distinct segments
  from $x$ to $y$, by Theorem \ref{unicita-giro} one of them, say
  $\ga_1 $, would have to pass through $0$.  But then $\ga_1$ would
  coincide with $ \ga_x * \ga_y^0$. This is absurd since $\ga_1$ is a
  segment.  This proves the first two assertions.  Since $\rcin\leq
  r_3 \leq r_2/2$, $x \in \Bd(0,r_2)$ and $d(x,y) < 2\rcin <
  r_2$. Thus the third assertion follows from the first and Lemma
  \ref{caratterizzo-c}.  The last fact follows from uniqueness by a
  standard use of Ascoli-Arzel\`a lemma.
%
%
\end{proof}

\begin{prop}
  \label{distanza-C1}
  The function $d(0, \cdot)$ is $\cl^1$ on $\Bd^*(0, \rcin)$. If
  $x\in \Bd^*(0, \rcin)$, the gradient of $d(0, \cdot )$ at $x$ is
  $\dga_x(d(0, x))$.
\end{prop}
\begin{proof}
  For $x$ and $\ga$ as above set $L=d(0,x)$ and $v=\dga(L)$. Choose
  $\eps$ such that $0< \eps < \min\{ L, \rcin - L\}$ and extend $\ga$
  to $[0,L+\eps]$ by setting $\ga(t)= \exp_x tv$ for $t\in (L,
  L+\eps]$.  Then $L+\eps \leq \rcin \leq r_2 $.  By
  Prop. \ref{inietto} the path $\ga\restr{[\delta, L+\eps]}$ is a
  segment for every $\delta>0$, hence $\ga$ is the unique segment from
  $0$ to $\ga(L+\eps)$. Set $x_1 = \ga(L-\eps)$ and $
  x_2=\ga(L+\eps)$.  Since $\eps = d(x_1,x) = d(x_2, x) < L=\inj{x}$
  the functions $d(x_1, \cdot)$ and $d(x_2, \cdot)$ are differentiable
  at $x$ with gradients $v$ and $-v$ respectively (see
  e.g. \cite[Prop. 4.8, p.108]{sakai-Riemannian}). So
  \begin{gather*}
    d(x_1, \exp_x w) = d(x_1, x) + g_x(v,w) + o(|w|)\\
    d(x_2, \exp_x w) = d(x_2, x) - g_x(v,w) + o(|w|).
  \end{gather*}
  By the triangle inequality
  \begin{gather*}
    \begin{aligned}
      d(0,  \exp_x w) &\leq d(0,x_1) + d(x_1,  \exp_x w) =\\
      &=d(0,x_1) +  d(x_1, x) + g_x(v,w) + o(|w|) = \\
      &=d(0,x) + g_x(v,w) + o(|w|)\\
      d(0, \exp_x w) &\geq d(0,x_2) - d(x_2, \exp_x w) =
      \\
      &= d(0,x_2) - \bigl [ d(x_2, x) - g_x(v,w) + o(|w|) \bigr] = \\
      &= d(0,x) + g_x(v,w) + o(|w|)
    \end{aligned}
    \\
    |d(0, \exp_x w) - d(0,x) - g_x(v,w)| = o (|w|).
  \end{gather*}
  This proves that $d(0,\cdot)$ is differentiable at $x$ with gradient
  $v$.  Next we show that the gradient is continuous. Indeed if
  $\{x_n\}$ is a sequence converging to $x\in \Xr$ and $\ga_n$ are
  segments from $0$ to $x_n$, then by Theorem \ref{W-Hoelder} and the
  theorem of Ascoli and Arzel\`a there is a subsequence $\ga_n^*$ that
  converges in the $\cl^1$-topology to the unique segment $\ga$ from
  $0$ to $x$. In particular $\dga_n^*(d(0,x_n)) \to
  \dga(d(0,x))$. Therefore the vector field $\nabla d(0,\cdot)$ is
  continuous on $\Bd^*(0,\rcin)$.
\end{proof}
We found the above argument for the differentiability of the distance
function in \cite[Prop. 6]{mccann-polar}.

\begin{lemma}
  \label{distanza-Riemanniana}
  Let $(M,g)$ be a Riemannian manifold with sectional curvature
  bounded above by $ \Kl\in \R$. Let $x$ and $y$ be points of $ M$
  that are connected by a unique segment $\ga(t) = \exp_x tv$, $v\in
  U_xM$ so that $\ga(t_0)=y$ and
  $$
  t_0=d(x,y) < \min\Bigl\{ c_x(v), \frac{ \pi}{2\sqrt{\Kl}}\Bigr\}
  $$
  where as usual $\sqrt{\Kl}=+\infty$ if $\Kl\leq 0$. Then the
  function $ d(x, \cdot)$ is smooth in a neighbourhood of $y$ and its
  Hessian at $y$ is positive semi-definite.
\end{lemma}
\begin{proof}
  This is a classical result in Riemannian geometry following from
  Rauch comparison theorem. It is commonly stated with stronger (and
  cleaner) hypotheses, but the usual proof, found e.g. in
  \cite[pp.151-153]{sakai-Riemannian} goes through without change with
  the above minimal assumptions. In fact $\exp_x$ is a diffeomorphism
  in a neighbourhood of $v\in T_xM$, so $d(x, \cdot)$ is smooth and one
  can compute its derivatives using Jacobi fields. The result then
  follows from Lemma 4.10 p.109 and Lemma 2.9 p.153 in
  \cite{sakai-Riemannian}, especially eq. (2.16) p.153. Notice that we
  are only interested in the first inequality in eq. (2.16) and this
  only depends on the upper bounds for the sectional curvature of $M$.
\end{proof}

\begin{prop}\label{medio-Whitehead}
  If $\alfa : [0, L] \ra \Bd^*(0, \rcin ) $ is a segment, the function
  $d(0, \alfa(\cdot))$ is convex on $[0,L]$.
\end{prop}
\begin{proof}
  Pick $s_0\in [0,L]$ and set $x=\alfa(s_0)$ and
  $x_n=\ga_x(1/n)$. Then $\uncobra(x_n) \equiv \uncobra(x)$ since
  $\ga_{x_n}$ is a piece of $\ga_x$. Since $\uncobra$ is continuous,
  there is an $\eps>0$ such that for any $s\in J: = (s_0 -\eps, s_0 +
  \eps) \cap [0,L]$
  \begin{equation*}
    \ango(\uncobra(x), \uncobra(\alfa(s)) ) < \frac{\pi}{m}.
  \end{equation*}
  By Prop.  \ref{piccolo-Whitehead} for any $s\in J$ there is a unique
  segment $\alfa_{n,s}: [0, d(x_n, \alfa(s))] \ra \Xr$, joining $x_n$
  to $\alfa(s)$, it is of the form $\alfa_{n,s}(t) = \exp_{x_n}
  tv_{n,s}$ and $d(x_n, \alfa(s)) < c_{x_n}(v_{n,s})$. So we can apply
  Lemma \ref{distanza-Riemanniana} to the effect that the function
  $u_n= d(x_n, \alfa(\cdot))$ is convex on $J$. Since $u_n \to d(0,
  \alfa(\cdot))$ uniformly, also the function $d(0, \alfa(\cdot))$ is
  convex on $J$. Since $t_0$ is arbitrary and convexity is a local
  condition, this proves convexity on the whole of $[0,L]$ as well.
\end{proof}
\begin{cor}
  \label{Whitehead}
  For any $r\in (0,\rcin)$ the ball $\Bd(0,r)$ is geodesically convex,
  that is: any segment whose endpoints lie in $ \Bd(0,r)$ is contained
  in $\Bd(0,r)$.
\end{cor}
\begin{proof}
  Let $\alfa :[0,L] \ra X$ be a segment with endpoints $x,y \in
  \Bd(0,r)$.  If $\alfa$ passes through the origin the assertion is
  obvious. Otherwise the function $u(t)= d(0, \alfa(t))$ is convex on
  $[0,L]$ by Proposition \ref{medio-Whitehead}. Since $x,y\in
  \Bd(0,r)$, $u(0)<r $ and $ u(L) <r$ so
  \begin{equation*}
    u(t) \leq \Bigl ( 1- \frac{t}{L} \Bigr ) u(0) + \frac{t}{L} u(L) < r.
  \end{equation*}
  Therefore $\alfa(t) \in \Bd(0,r)$ for any $t\in [0,L]$.
\end{proof}
Now choose a number $\rsei \in (0, \rcin)$ and set
\begin{equation}
  \label{eq:def-C}
  C=\{x\in X: d(0,x)=\rsei\}.
\end{equation}
It follows from Proposition \ref {distanza-C1} that $C$ is a smooth
1-dimensional submanifold of $\Xr$. Since it is compact, it is
diffeomorphic to $S^1$.  The interior of $C$ is $\Bd(0,\rsei)$ which
is thus a topological disc. Let $\sigma: \R \ra C$ be a positively
oriented $\cl^1$ periodic parametrisation of $C$ of period $1$. Since
$\sigma$ is positively oriented the vector $J\dot{\sigma}$ points
inside $\Bd(0,\rsei)$.
\begin{lemma}\label{nonconsto}
  The maps $\cobra\circ \sigma$ and $\uncobra\circ \sigma$ are not
  constant.
\end{lemma}
\begin{proof}
  Since $ \cobra(x) = (u(\uncobra(x)), 0, \ldots{}, 0) $ it is enough
  to prove that $\cobra\circ\sigma $ is not constant. Assume by
  contradiction that $\cobra(x) \equiv v$ for any $x\in C$.  Since the
  range of $\cobra$ is contained in $C_0(X)$,
  $\pi\circ\cobra=\cobra$. Using \eqref{eq:angoletto-5} we get for $
  x\in C$, $ s\in [0, \rsei]$
  \begin{gather*}
    \ango\bigl(v, \pi(\dga_x(s)\bigr) \leq \ango( v, \cobra(x)) +
    \ango( \pi(\dga_x(0)), \pi(\dga_x(s))) \leq
    \deltone \\
    \pi(\dga_x(s))\cdot v > 0 \\
    \pi(x) \cdot v = \pi(\ga_x(\rsei)) \cdot v = \int_0^{\rsei}
    \pi(\dga_x(s)) \cdot v \, ds > 0.
  \end{gather*}
  Therefore $\pi(C)=\pi\sigma([0,1])$ would be contained in the
  half-plane $\{z\in \C: z\cdot v > 0\}$ and $\pi\circ \sigma
  \restr{[0,1]}$ would be null-homotopic in $\Delta^*$. This is
  impossible since $\sigma\restr{[0,1]}$ generates $\pi_1(\Xr,
  \sigma(0))$ and $\pi: \Xr \ra \Delta^*$ is an $m:1$ covering.
\end{proof}

\begin{definiz}
  For $t_0, t_1 \in \R$ set
  \begin{gather}
    \nonumber
    \intT=\{ s e^{2\pi it} \in \C: s \in (0,1) , t\in (t_0,t_1)\} \\
    \barT=\{ s e^{2\pi it} \in \C: s \in [0,1] , t\in [t_0,t_1]\}
    \nonumber{}
    \\
    \label{eq:def-marmo}
    \marmo : \barT \ra X \qquad \marmo (se^{it}) = \camel (\sigma(t),
    s)
    \\
    \label{eq:def-sector}
    \sect(t_0,t_1) = \marmo (\intT) \qquad \sect[t_0,t_1] = \marmo
    (\barT).
  \end{gather}
\end{definiz}
Since $d(\sigma(t),0)=\rsei$ for any $t\in \R$ 
\begin{gather*}
  \marmo (se^{it}) 
  = \ga_{\sigma(t)} (\rsei s).
\end{gather*}
\begin{lemma}\label{dentro-dentro}
  If $t_0< t_1 < t_0+1$ the map $\marmo$ is a homeomorphism of $\barT$
  onto $\sect[t_0,t_1]$, $ \inte \sect[t_0,t_1] = \sect(t_0,t_1) $ and
  \begin{gather}
    \label{eq:bordo-marmo}
    \partial \sect[t_0,t_1] = \sigma([t_0,t_1]) \cup \Im
    \ga_{\sigma(t_0)}\cup \Im \ga_{\sigma(t_1)}.
  \end{gather}
\end{lemma}
\begin{proof}
  Continuity of $\marmo$ follows from Proposition \ref{F-cont}. We
  prove that it is injective. Let $s e^{2\pi it}, s' e^{2\pi it'} \in
  T$, be such that $\marmo(s^{2\pi it}) =\marmo (s'e^{2\pi it'})=y$.
  If $s = 0$ then
  \begin{gather*}
    y= \ga_{\sigma(t') } (\rsei s' )
    = 0
  \end{gather*}
  so $ s' = 0$ as well and $se^{2\pi it} = s' e^{2\pi it'}=0$.
  If $s, s' >0$, write $x=\sigma(t)$, $x'=\sigma(t')$. Then
  \begin{equation*}
    \ga_{x} (\rsei s) = 
    \ga_{x'}
    ( \rsei s') = y.
  \end{equation*}
  So $\rsei s=d(0,y)= \rsei s'$ and $s=s'$.  Moreover, from Theorem
  \ref{unicita}, we get $\ga_{x} (t) =\ga_{x'}(t)$ for $t\in
  [0,d(y,0)]$ and by the unique continuation of geodesics also for
  $t\in [d(y,0),\rsei]$. Hence $x=\ga_x(\rsei)=\ga_{x'}(\rsei)=x'$ and
  $t=t'$. This shows that $\marmo$ is injective and therefore a
  homeomorphism of $T$ onto its image $\marmo(T)$. Since $T$ is
  homeomorphic to a closed disk, Brouwer theorem on the invariance of
  the domain and of the boundary (see e.g.
  \cite[p.205f]{prasolov-homology}) implies that $\inte
  \marmo(T)=\marmo (\inte T)$ and $\partial \marmo (T) = \marmo(T) -
  \marmo (\inte T) = \marmo (\partial T) = \sigma([t_0,t_1]) \cup \Im
  \ga_{\sigma(t_0)}\cup \Im \ga_{\sigma(t_1)}$.
\end{proof}
Set $ p(t)=e^{2\pi i t}$ and let $\solle$ a lifting of $\uncobra\circ
\sigma$:
\begin{gather}
  \label{eq:def-solle}
\begin{diagram}
  \node{ \R} \arrow[2]{e,t}{\solle} \arrow{s,l}{\sigma}
  \arrow{ese,t}{\uncobra\circ \sigma} \node[2]{\R} \arrow{s,r}{p}
  \\
  \node{ C } \arrow[2]{e,b}{\uncobra} \node[2]{ S^1}
\end{diagram}
\end{gather}

\begin{lemma}\label{surietto}
  The function $\solle$ is monotone increasing and $\uncobra(C)=S^1$.
\end{lemma}
\begin{proof}
  Mark that we are not saying that $\solle$ is \enf{strictly}
  increasing.  Let $t_0,t_1\in \R$ be such that $t_0< t_1 < t_0 +1
  $. Set $x_0=\sigma(t_0)$, $x_1=\sigma(t_1)$,
  $R=\phi\meno\bigl(\sect(t_0,t_1)\bigr)$, see \eqref{eq:def-marmo}.
  Since $\phi\meno$ is an orientation preserving homeomorphism of
  class $\cl^1$ outside the origin, it follows from Lemma
  \ref{dentro-dentro} and Prop. \ref{rego-in-D} that $R$ is a region
  of $\Delta$ homeomorphic to a disk with piecewise $\cl^1$ boundary
  \begin{gather*}
    \partial R = \phi\meno \sigma ([t_0, t_1]) \cup \Im
    \unga_{\sigma(t_0)}\cup \Im\unga_{\sigma(t_1)}.
  \end{gather*}
  Moreover $R$ lies on the left of $\unga_{\sigma(t_0)}$ and on the
  right of $\unga_{\sigma(t_1)}$ and for $t\in (t_0, t_1)$ the path
  $\unga_{\sigma(t)}$ lies inside $R$.  Accordingly its tangent vector
  \begin{gather*}
    \dunga_{\sigma(t)}(0) = e^{2\pi i \solle(t)}
  \end{gather*}
  points inside $R$.  Let $\psi\in[0,1)$ be such that $\uncobra(x_1) =
  e^{2\pi i\psi} \uncobra(x_0)$.
  Since $\unga_{\sigma(t_0)}(0)=\uncobra(x_0)$ and $
  \unga_{\sigma(t_1)}(0)=\uncobra(x_1)$ the unit tangent vectors at
  $0$ pointing inside $R$ are exactly those of the form
  $e^{i\theta}\unga_{x_0}(0)$ with $ \theta \in [0,2\pi \psi]$.  So
  $$
  \solle [t_0,t_1] \subset \bigsqcup_{k\in \Zeta} [\solle(t_0) + k,
  \solle(t_0) + \psi +k ].
  $$
  Since $\solle$ is continuous, we have $\solle[t_0,t_1 ] =
  [\solle(t_0) , \solle(t_0) + \psi]$ and $\solle(t_1) = \solle(t_0) +
  \psi$.  This proves that $\solle(t_1) \geq \solle(t_0) $. It follows
  that $\solle $ is increasing on the real line.  Since $p\solle (1) =
  p \solle (0)$, there is $k\in \Zeta$ such that $\solle (1) = \solle
  (0) + k$.  By the uniqueness of the lifting $\solle(t+1)=\solle(t) +
  k$ for any $t\in \R$.  $k \geq 0$ because $\solle$ is increasing.
  If $k=0$ then $\solle$ would be constant on $ [0,1]$ and so on the
  whole real line. But this is not the case by Lemma \ref
  {nonconsto}. Therefore $k>0$.  It follows that $\uncobra $ is
  surjective.
\end{proof}

\begin{lemma} 
  \label{giro-poco-ma-giro}
  If $x\in \Bd^*(0,\rsei)$ then
  \begin{gather*}
    \ango(\pi(x), \cobra(x)) < \deltone.
  \end{gather*}
\end{lemma}
\begin{proof}
  For $x\in \Bd^*(0,\rsei)$ set $L=d(0,x)$ and let $\ga_x$ be as in
  \eqref{eq:def-camel}. The set $E=\{w\in \C^*: \ango\bigl
  (w,\cobra(x) \bigr) < \pi/8
  \}$ is a convex cone. Since $\cobra(x)=\dga_x(0) = \pi(\dga_x(0))$,
  it follows from \eqref{eq:angoletto-5} that $\pi(\dga_x(s))\in E$
  for any $s\in [0,L]$. Thus
  \begin{gather*}
    \pi(x) = \int_0^L \pi(\dga_x(s)) ds \in E.
  \end{gather*}
\end{proof}

\begin{teo} \label{settore-convesso} Let $t_0 , t_1 \in \R$ be such
  that $t_0< t_1$ and
  \begin{gather*}
    \solle(t_1) < \solle(t_0) + \frac{1}{2m}.
  \end{gather*}
  Then for any $x,y\in \sect[t_0,t_1]$ there is a unique segment
  joining $x$ to $y$ and it is contained in
  $\sect[t_0,t_1]$. Therefore $\sect(t_0,t_1)$ and $\sect[t_0,t_1]$
  are geodesically convex subsets of $(X,d)$.
\end{teo}
\begin{figure}[h]
  \begin{center}
    \includegraphics{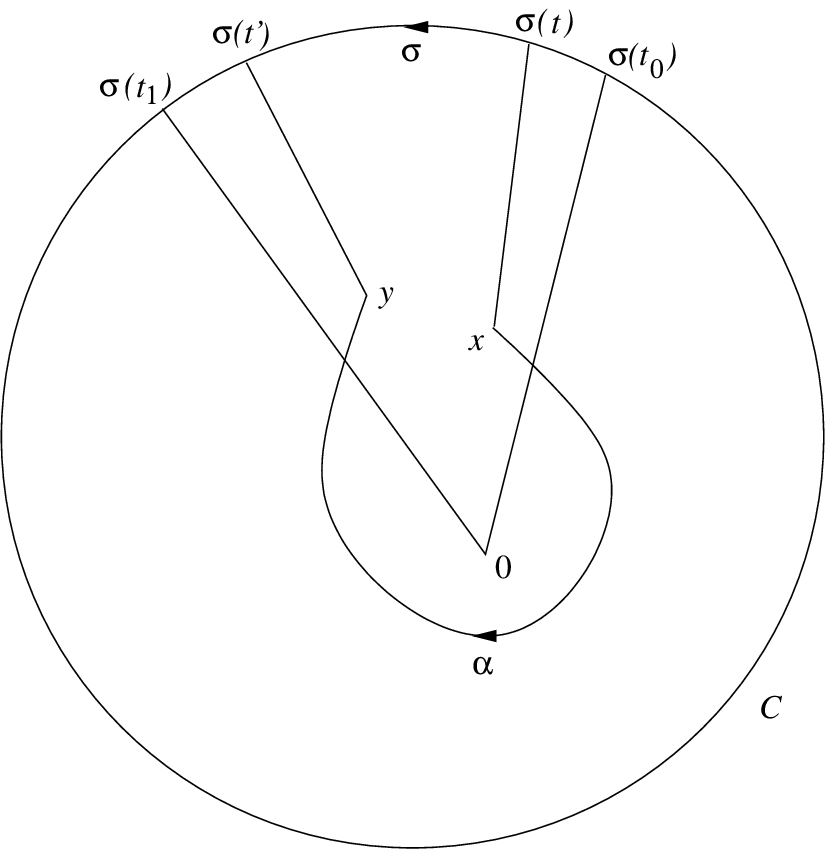}
  \end {center}
  \caption{}
  \label{torta}
\end{figure}
\begin{proof}
  Since $\sect(t_0,t_1)$ can be exhausted by sets of the form
  $\sect[t_0+\delta, t_1 - \delta]$ with $(t_1 - \delta) - (t_0 +
  \delta) < \um$ it is enough to prove the convexity of $\sect[t_0,
  t_1]$.  Let $x$ and $ y$ be points in $\sect[t_0, t_1]$.  If either
  $x=0$ or $y=0$ the claim is immediate from the definition of
  $\sect[t_0,t_1]$. Otherwise we can assume that
  \begin{gather*}
    \begin{gathered}
      x=\camel(\sigma(t),s) = \ga_{\sigma(t)}(\rsei s )
      \\
      y=\camel ( \sigma(t'), s') = \ga_{\sigma(t')}(\rsei s' )
    \end{gathered}
    \qquad
    \begin{gathered}
      t_0 \leq t \leq t' \leq t_1 \\
      s,s'\in (0,1].
    \end{gathered}
  \end{gather*}
  Since $\solle$ is monotone $ \solle (t_0) \leq \solle(t) \leq
  \solle(t') \leq \solle (t_1) $ and
  \begin{gather*}
    \ango\bigl(\uncobra(x), \uncobra(y)\bigr) = 2\pi | \solle(t) -
    \solle(t')| \leq 2\pi \bigl ( \solle(t_1) - \solle(t_0) \bigr ) <
    \frac{\pi}{m}.
  \end{gather*}
  It follows from Lemma \ref{ball-hopf} and
  Prop. \ref{piccolo-Whitehead} that there is a unique segment $\alfa:
  [0,L] \ra X$ joining $x$ to $y$ (so $L=d(x,y)$).  We need to prove
  that $\alfa ( [0,L]) \subset \sect[t_0,t_1]$.  Assume that it is
  not. Then $\alfa$ has to cross $\partial \sect[t_0,t_1]$ at least
  twice.  Since $d(0,x) < \rsei$ and $ d(0, y) < \rsei$, we have
  $\alfa([0,L]) \subset \Bd(0,\rsei)$ by Corollary \ref{Whitehead}.
  It follows from \eqref{eq:bordo-marmo} that the set $ \Im \alfa\cap
  \bigl (\Im \ga_{\sigma(t_0)} \cup \Im \ga_{\sigma(t_1)}\bigr ) $
  contains at least two points.  On the other hand $\alfa$ cannot
  cross the path $\ga_{\sigma(t_i)}$ more than once:
  otherwise by Theorem \ref{unicita-giro} it would coincide with some
  prolongation of $\ga_{\sigma(t_i)}$. This proves that $\alfa$
  crosses each of the paths $\ga_{\sigma(t_i)}$ exactly once. Define
  $\beta: [0,3]\ra X$ by
  \begin{gather*}
    \beta (\tau) =
    \begin{cases}
      \ga_{\sigma(t)}\bigl( \rsei(s + \tau (1-s))\bigr )
      & \text{if } \tau \in [0,1]\\
      \sigma (t + (\tau -1)(t'-t) & \text{if }\tau \in [1,2]\\
      \ga_{\sigma(t')}\bigl ( \rsei (1 + (\tau -2 ) (s'-1))\bigr ) &
      \text{if } \tau \in [2,3].
    \end{cases}
  \end{gather*}
  Then $\zeta = \beta * \alfa^0$ is a simple closed curve, $[\zeta]$
  is the positive generator of $\pi_1(\Xr, x)$ and $W(\pi\circ
  \zeta)=m$.  By Lemma \ref{non-giro} $W(\pi\circ \alfa) < 1$, so
  \begin{gather}
    \label{eq:W-contraddicendo}
    W(\pi\circ \beta) > m-1 \geq 1.
  \end{gather}
  Let $ \tau' : [0,3] \ra [t,t']$ be the function
  \begin{gather*}
    \tau'(\tau)=\begin{cases}
      t & \text{if } \tau \in [0,1]\\
      t+ (\tau -1)(t'-t) & \text{if } \tau \in [1,2]\\
      t' & \text{if } \tau \in [2,3].
    \end{cases}
  \end{gather*}
  Clearly $\cobra\bigl(\beta(\tau)\bigr) = \cobra\bigl(\sigma(\tau'(\tau))\bigr)
  $.  By Lemma \ref{giro-poco-ma-giro}
  \begin{gather*}
    \ango\bigl (\pi\circ\beta(\tau) , \cobra(\sigma(\tau'(\tau)) \bigr) <
    \deltone
  \end{gather*}
  for any $\tau\in[0,3]$.  Let $t_2\in [t,t']$ be such that
  \begin{gather*}
    \solle(t_2) = \frac{\solle(t) + \solle(t')}{2}
  \end{gather*}
  and set $w=\cobra(\sigma(t_2))$.  Then for any $\tau'\in [t,t']$
  \begin{gather*}
    |\solle(\tau') - \solle(t_2) | \leq \frac{|\solle(t') -
      \solle(t)|}{2} \leq \frac{\solle(t_1) - \solle(t_0) }{2} <
    \frac{1}{4m}
  \end{gather*}
  so for $\tau \in [0,3]$
  \begin{gather*}
    \ango\bigl(\pi\circ\beta(\tau), w\bigr) \leq \ango\bigl(\pi\circ\beta(\tau),
    \cobra(\tau'(\tau))\bigr) +
    \ango\bigl(\cobra(\tau'(\tau)),w\bigr) < \\
    < \deltone + 2\pi m |\solle(\tau') - \solle (t_2) | <
    \frac{5\pi}{8}.
  \end{gather*}
  This shows that $W(\pi\circ \beta) \leq 1$, contradicting
  \eqref{eq:W-contraddicendo}. Therefore $\alfa([0,L]) \subset
  \Bd(0,\rsei)$ as claimed.
\end{proof}

\section{Alexandrov curvature}

\label{Alex-sec}

In this section we will finally conclude the proof of Theorem
\ref{main}.  We start by recalling the basic definitions related to
upper curvature bounds for a metric space in the sense of
A.D.Alexandrov.  Next we will come back to the setting considered in
\Sec\Sec \ref{sec-regu}--\ref{convex-section} and we will prove that
$\Bd(0,\rsei)$ is a \cat--space
(Thm. \ref{caso-irriducibile}). Theorem \ref{main} follows almost
immediately from this.

\bigskip{}

A thorough treatment of the intrinsic geometry of metric spaces, and
especially of curvature bounds in the sense of Alexandr Danilovich
Alexandrov can be found in the books \cite{aleksandrov-zalgaller},
\cite{rinow}, \cite{aleksandrov-berestovskij-nikolaev-eng},
\cite{berestovskij-nikolaev-enci}, \cite{ballmann-nonp},
\cite{bridson-haefliger-libro}, \cite{burago-burago-ivanov}.  We
mostly follow 
\cite{bridson-haefliger-libro}.

Let $(X,d)$ denote an arbitrary metric space with intrinsic metric.
Given two segments $\alfa$ and $\beta$ in $X$ with
$\alfa(0)=\beta(0)=x$ the \enf{Alexandrov (upper)
  angle} 
is defined as
\begin{gather}
  \label{eq:def-angle}
  \alex_x(\alfa, \beta) = \limsup_{t, t' \to 0} \, \arccos \ 
  \frac{ t^2 + (t')^2 - d(\alfa(t), \beta(t'))^2 }{2tt'} 
  .
\end{gather}
Fix $\Kl\in \R$. Set $D_\Kl= +\infty$ if $\Kl \leq 0$ and $D_k =
\pi/\sqrt{\Kl}$ otherwise.
Let $M^2_\Kl$ denote the complete Riemannian surface
with constant curvature $\Kl$.  A \enf{triangle} $T=\T(xyz)$ in $X$ is
a triple of points $x,y,z$ together with a choice of three segments
connecting them. A \enf{comparison triangle}
is a triangle $\overline{T} =\T( \bar{x} \bar{y} \bar{z})$ in
$M^2_\Kl$ such that corresponding edges have equal length. We will
occasionally let $T$ denote also the union of the edges.

  
\begin{definiz} \label{angularis-condicio} We say that the \enf{angle
    condition} holds for a triangle $T$ in a metric space, if the
  Alexandrov angle between any two edges of $T$ is less or equal than
  the angle at the corresponding vertex in a comparison triangle
  $\overline{T}$ in $M^2_\Kl$.  A metric space $(X,d)$ is called
  \enf{\cat--space} if (1) the metric is intrinsic, (2) any pair of
  points $x,y\in X$ with $d(x,y) < 2D_\Kl$ is connected by a segment
  and (3) the angle condition holds for any triangle in $X$.
  A metric space has \enf{ curvature $\leq \Kl$ (in the sense of
    Alexandrov)} if 
  for every $x\in X$ there is $r_x>0$ such that the ball of centre $x$
  and radius $r_x$
  endowed with the induced metric is a \cat--space.
\end{definiz}

\begin{prop}\label{cara-cat}
  Let $\Kl\in \R$ and let $(X,d)$ be a $D_\Kl$--geodesic metric space
  (this means that any pair of points a distance less than $D_\Kl$
  apart are connected by a segment). Then $(X,d)$ is a \cat--space if
  and only if
  for any triangle $T$ in $X$ with perimeter less than $2D_\Kl$ the
  following condition holds: for $x, y \in T$ let $\bar{x}$ and
  $\bar{y}$ denote the corresponding points on a comparison triangle
  in $M^2_\Kl$; then $d(x,y) \leq d(\bar{x},\bar{y})$.
\end{prop}
See \cite[p.161]{bridson-haefliger-libro}.

\begin{prop}\label{AlexRiem}
  Let $(M,g)$ be a Riemannian manifold with sectional curvature
  bounded above by $\Kl$. Then $M$ provided with the Riemannian
  distance is a metric space of curvature $\leq \Kl$ in the sense of
  Alexandrov.
\end{prop}
For a proof see e.g. \cite[Thm. 2.7.6 p. 219]{klingenberg-libro} or
\cite[Thm. 1A.6 p.173]{bridson-haefliger-libro}).

\begin{prop}\label{angolo-piatto}
  Let $(X,d)$ be a \cat--space and let $\alfa: [0,a]\ra X$ and $\beta
  : [0,b]\ra$ be two segments with $\alfa(0)=\beta(0)=x$ and
  $\alex_x(\alfa, \beta) = \pi$. Then $\alfa^0 * \beta$ is a segment.
\end{prop}
See \cite[Prop. 9.1.17(4) p.313]{burago-burago-ivanov}.

\begin{lemma}
  \label{gluing-triangles}
  Let $(X,d)$ be a $D_\Kl$--geodesic space and let $T=\T(xy_1y_2)$ be
  triangle with perimeter $< 2D_\Kl$ and distinct vertices. Fix a
  point $z$ on the segment from $y_1$ to $y_2$ and a segment from $x $
  to $z$.  In this way we get two triangles $T_1=\T(xy_1z)$ and
  $T_2=\T(xzy_2)$ with a common edge.  If the angle condition holds
  for both $T_1$ and $T_2$ then it also holds for $T$. 
\end{lemma}
This is the gluing lemma of \cite[p.199]{bridson-haefliger-libro}.

\begin{prop}\label{cat-in-Alex} Let $(X,d)$ be a metric space of
  curvature $\leq \Kl$.  Assume that for every pair of points $x,y \in
  X$ with $d(x,y) < D_\Kl$ there is a unique segment $\alfa_{x,y}$
  which depends continuously on $(x,y)$. Then $X$ is a \cat--space.
\end{prop}
See \cite[Prop. 4.9 p.199]{bridson-haefliger-libro}.

\begin{prop}
  Let $(X, d)$ be a \cat--space and let $\alfa$ and $\beta$ be
  segments with $\alfa(0)=\beta(0)=x$. Then the $\limsup$ in
  \eqref{eq:def-angle} is in fact a limit. Therefore
  \begin{gather}
    \label{eq:angolo-forte}
    \alex_x(\alfa, \beta) = \lim_{t\to 0} \arccos \frac{2t^2 -
      d(\alfa(t), \beta(t))^2 }{2t^2} = 2 \lim_{t\to 0} \arcsin \frac{
      d(\alfa(t), \beta(t)) }{2t} .
  \end{gather}
\end{prop}
If $(X,d)$ is a uniquely geodesic metric space we denote by $[x,y]$
the segment from $x$ to $y$ and by $\angle_x(y,z)$ the Alexandrov
angle between the segments $[x,y]$ and $[x,z]$.
\begin{prop}\label{angolo-continuo}
  If $(X,d)$ is a \cat--space the function $(x,y,z) \mapsto
  \angle_x(y,z)$ is upper semicontinuous on the set of triples $(x,y,z)$ with
  $d(x,y), d(x,z) < D_\Kl$. For fixed $x$ the function $(y,z) \mapsto
  \angle_x(y,z)$ is continuous.
\end{prop}
See \cite[pp.184-185]{bridson-haefliger-libro}.

\bigskip{}

Let us now come back to the setting and the notation of
\Sec\Sec \ref{sec-regu}--\ref{convex-section}.  For $t\in \R$ set
\begin{gather*}
  \opt= \sup \Bigl \{ \tau > t : \solle(\tau) < \solle(t) +
  \frac{1}{2m} \Bigr \} \\
  \opmt= \inf \Bigl \{ \tau < t : \solle(\tau) > \solle(t) -
  \frac{1}{2m} \Bigr \}.
\end{gather*}
Clearly
\begin{gather}
  \label{eq:sollepm}
  \solle(\oppmt) = \solle (t) \pm \frac{1}{2m}.
\end{gather}
From the monotonicity of $\solle$ it follows that if $t''<t< t'$ then
\begin{gather}
  \label{eq:ttprimo}
  \begin{gathered}
    \solle(t') < \solle(t) + \frac{1}{2m} \Longleftrightarrow t' < \opt \\
    \solle(t'') > \solle(t) - \frac{1}{2m} \Longleftrightarrow t'' >
    \opmt.
  \end{gathered}
\end{gather}
\begin{prop}
  \label{settore-convesso-massimale}
  For any $t\in \R$ the sectors $\sect(t,\opt)$ and $\sect(\opmt, t)$
  are geodesically convex subsets of $(X,d)$. Moreover
  $\sect(t,\opt)$, $\sect(\opmt,t)$, $\sect[t,\opt]$ and
  $\sect[\opmt,t]$ provided with the distance induced from $(X,d)$ are
  \cat--spaces.
\end{prop}
\begin{proof}
  We consider only $\sect(t,\opt)$ and $\sect[t,\opt]$.  If $x_0,x_1
  \in \sect(t,\opt)$ then by \eqref{eq:def-sector} $x_i= \camel
  (\sigma(t_i),s_i)$ with $t_i \in (t, \opt) $.  Assume $t_0 < t_1$
  and set
  \begin{gather*}
    t_0' = \frac{t_0+t}{2} \qquad t_1' = \frac{t_1+\opt}{2}.
  \end{gather*}
  Then $t< t_0' < t_0 < t_1 < t_1' < \opt$ and
  \begin{gather*}
    \solle(t_1') < \solle(t) + \frac{1}{2m} \leq \solle(t_0') +
    \frac{1}{2m}.
  \end{gather*}
  By Theorem \ref{settore-convesso} there is a unique segment from
  $x_0$ to $x_1$, and it is contained in $\sect(t_0', t_1') \subset
  \sect(t,\opt)$. It follows that $\sect(t,\opt)$ is a geodesically
  convex subset of $(X,d)$.  In particular $\sect(t,\opt)$ provided
  with the distance induced from $(X,d)$ is a geodesic metric space. By
  continuity the same holds for $\sect[t,\opt] = \overline{\sect(t,
    \opt)}$: for any $x_0, x_1 \in \sect[t,\opt]$ there is \enf{at
    least one} segment from $x_0$ to $x_1$ that is contained in
  $\sect[t,\opt]$.  Moreover the induced distance on $\sect(t, \opt)$
  coincides with the Riemannian distance of the smooth Riemannian
  surface $(\sect(t, \opt), g\restr{\sect(t,\opt)})$, whose Gaussian
  curvature is everywhere $\leq \Kl$.
  Prop. \ref{AlexRiem} ensures that $(\sect(t,\opt), d)$ is a metric
  space with curvature $\leq \Kl$ in the sense of Alexandrov. 
   By Thm.  \ref{settore-convesso}
  it is uniquely geodesic and by Prop. \ref{piccolo-Whitehead}
  segments in $\sect(t,\opt)$ depend continuously on the
  endpoints. Thus Prop. \ref{cat-in-Alex} ensures that
  $(\sect(t,\opt), d)$ is a \cat--space. Since any triangle in  $\sect[t,\opt]$
  is a limit of triangles in $\sect(t, \opt)$ a continuity argument
  applied to the condition in Prop. \ref{cara-cat} yields
that  $\sect[t,\opt]$ is a
  \cat--space too.
\end{proof}

\begin{prop} \label{Alex-angle} If $\solle(t_1) < \solle(t_0) +
  \frac{1}{2m}$ then
  \begin{gather}
    \label{eq:Alex-angle}
    \alex_0\bigl (\ga_{\sigma(t_0)}, \ga_{\sigma(t_1)} \bigr) = 2\pi m
    \bigl (\solle(t_1) - \solle(t_0) \bigr).
  \end{gather}
\end{prop}
\begin{proof}
  Both $\ga_{\sigma(t_0)}$ and $ \ga_{\sigma(t_1)}$ are segments
  contained in the \cat--space $\sect[t_0,t_0^+]$. By
  \eqref{eq:angolo-forte} and \eqref{eq:meno-di-1} their Alexandrov
  angle is
  \begin{gather*}
    \alex_0\bigl (\ga_{\sigma(t_0)}, \ga_{\sigma(t_1) } \bigr ) = 2
    \lim_{t\to 0} \arcsin \frac{ d(\ga_{\sigma(t_0)}(t),
      \ga_{\sigma(t_1)}(t)) }{2t} = \\
    =\ango \bigr (\dga_{\sigma(t_0)}(0), \dga_{\sigma(t_1)}(0) \bigl )
    = \ango (\cobra(\sigma(t_0)), \cobra(\sigma(t_1))) .
  \end{gather*}
  Since $\uncobra (\sigma(t_i)) = e^{2\pi i \solle (t_i)}$, $\cobra
  (\sigma(t_i)) = (e^{2\pi m i \solle(t_i)}, 0, \ldots{}, 0)$ 
and $2\pi m |\solle(t_0) - \solle(t_1) | < \pi$ 
we get
  \begin{gather*}
    \ango (\cobra(\sigma(t_0)), \cobra(\sigma(t_1))) = 2 \pi m \big
    |\solle(t_0) - \solle(t_1) \big|.
  \end{gather*}
\end{proof}

\begin{prop}\label{segmenti-dritti}
  For any $t\in \R$ both $\ga_{\sigma(t)}* \ga_{\sigma(\opt)}^0$ and
  $\ga_{\sigma(t)}* \ga_{\sigma(\opmt)}^0$ are shortest paths.
\end{prop}
\begin{proof}
  Consider the first path. Thanks to
  Props. \ref{settore-convesso-massimale} and \ref{angolo-piatto} it
  is enough to show that $\alex_0\bigl (\ga_{\sigma(t)},
  \ga_{\sigma(\opt)} \bigr ) = \pi$. Indeed by
  Props. \ref{angolo-continuo} and \ref{Alex-angle} and
  \eqref{eq:sollepm}
  \begin{gather*}
    \alex_0\bigl (\ga_{\sigma(t)}, \ga_{\sigma(\opt)} \bigr ) = \lim
    _{\tau < \opt, \ \tau \to \opt} \alex_0\bigl (\ga_{\sigma(t)},
    \ga_{\sigma(\tau)} \bigr ) = \\
    =\lim_{\tau < \opt, \ \tau \to \opt} 2\pi m \bigl (\solle(\tau) - \solle(t)
    \bigr) = 2\pi m \bigl (\solle(\opt) - \solle(t) \bigr) = \pi.
  \end{gather*}
\end{proof}

Now we are able to control uniqueness of geodesics in general.

\begin{teo}\label{unicita-superfinale}
  For any $x,y\in \Bd(0,\rsei)$ there is a unique segment $\alfa_{x,y}$
  from $x$ to $y$ and it depends continuously on its endpoints. If 
$x\neq 0 $ and $y\neq 0$
the segment $\alfa_{x,y}$ passes
  through $0$ if and only if
  \begin{gather}
    \label{eq:evito}
    \ango \bigl (\uncobra(x), \uncobra(y) \bigr ) \geq \frac{\pi}{m}.
  \end{gather}
\end{teo}
\begin{proof}
  If one of the points is $0$ uniqueness is proved in Theorem
  \ref{unicita}.  If $x,y$ are two distinct points in $ \Bd^*(0,\rsei)$
  by interchanging them if necessary we can write them as
  \begin{gather*}
    x=\camel (\sigma(t), s) \qquad y=\camel (\sigma(t'), s')
  \end{gather*}
  with $0\leq t \leq t' <1$.
  We distinguish three cases according to the position
  of $t'$ with respect to $\opmt, t, \opt$. \\
  \noindent{} 1.  Assume first that $t' \in (t, \opt)$.  Then
  \begin{gather*}
    \solle(t') < \solle (t) + \frac{1}{2m} .
  \end{gather*}
\begin{figure}[h]
  \begin{center}
    \includegraphics{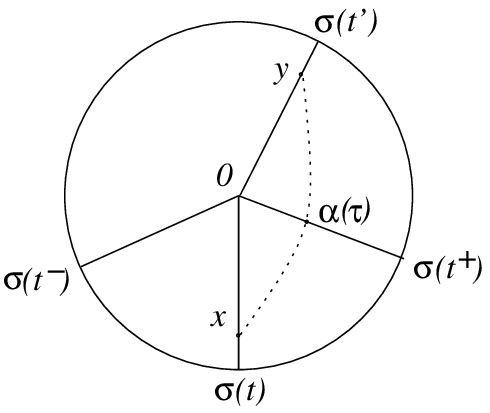}
  \end {center}
  \caption{}
  \label{sentieri}
\end{figure}
  It follows from Theorem \ref{settore-convesso} that the segment is
  unique, is contained in $\sect[t,t']\subset \sect[t, \opt]$ and does
  not pass through $0$ because of Lemma \ref{geodetiche-angolose}.
  \\
  \noindent{} 2.  If $t' > 1+\opmt$ then 
  \begin{gather*}
    \solle(t) < \solle (t'-1) + \frac{1}{2m} .
  \end{gather*}
  Since $y=\camel (\sigma(t'-1), s')$ the same argument proves that
  the segment is unique, is contained in $\sect[\opmt, t]$ and does
  not pass through $0$.\\
  \noindent{}
  3.  Finally assume that $ t' \in [\opt, 1+\opmt]$. This condition is
  just a restatement of \eqref{eq:evito}. In this case it is enough to
  prove that any segment $\alfa$ from $x$ to $y$ necessarily passes
  through the origin: Theorem \ref{unicita} then yields $\alfa = \ga_x
  ^0* \ga_y$ and in particular uniqueness.  Assume by contradiction
  that $\alfa $ does not pass through $0$ (see Fig. \ref{sentieri}).
  Then it has to cross either $\ga_{\sigma(\opt)}$ or
  $\ga_{\sigma(\opmt)}$.  Assume for example that $\alfa(\tau) =
  \ga_{\sigma(\opt)}(\tau')$
for some $\tau \in [0,d(x,y)]$,  $\tau'\in [0, \rsei]$. 
By 
Prop. \ref{segmenti-dritti}
  $\ga_x^0 *\bigl( \ga_{\sigma(\opt)} \restr{[0,\tau']} \bigr )$ is a segment.
So $\alfa\restr{[0,\tau]}$ and
  $\ga_x^0 * \bigl ( \ga_{\sigma(\opt)} \restr{[0,\tau']} \bigr)$ are two segments contained in
  $\sect[t,\opt]$ with same endpoints. By
  Prop. \ref{settore-convesso-massimale} $\sect[t,\opt]$ is a
  \cat--space, hence uniquely geodesic. Therefore
  $\alfa\restr{[0,\tau]}$ and $\ga_x^0 *  \bigl (\ga_{\sigma(\opt)}
  \restr{[0,\tau']}\bigr )$ must coincide, contrary to the assumption that
  $\alfa$ does not pass through the origin.

  Continuous dependence from the endpoints follows from uniqueness by
  Ascoli-Arzel\`a theorem.
\end{proof}

\begin{teo}
  \label{caso-irriducibile}
  The ball $\Bd(0, \rsei)$ provided with the distance induced from
  $(X,d)$ is a \cat--space.
\end{teo}
\begin{proof}
  We will show that any geodesic triangle $T=\T(xyz)$ contained in
  $\Bd(0,\rsei)$ satisfies the angle condition,
  Def. \ref{angularis-condicio}. We
  distinguish various cases. \\
  \noindent{}
  1.  Suppose first that the origin is a vertex, say $z=0$. If $ \ango
  \bigl (\uncobra(x), \uncobra(y) \bigr) < \pi/m$
  by interchanging if necessary $x$ and $y$ we can assume that
  $x=\camel (\sigma(t), s)$ and $y=\camel (\sigma(t'), s')$ with $t
  \leq t' < \opt$.  Then $T \subset \sect[t,t']$.  Since $
  \sect[t,t']$ is a \cat--space, the angle
  condition holds for $T$.\\
  \noindent{}
  2.  If $z=0$ and $ \ango \bigl (\uncobra(x), \uncobra(y) \bigr )
  \geq \pi/m $, then $\alfa_{x,y}=\ga_x^0 * \ga_y$ by Theorem
  \ref{unicita-superfinale}. So $T$ is
  degenerate and trivially satisfies the angle condition.\\
  \noindent{}
  3. Next assume that $0$ belongs to some edge but is not a vertex.
  Say $0$ lies on the edge $[x, y]$. By the above both triangles
  $\T(0xz)$ and $\T(0yz)$ satisfy the angle condition. By  Lemma
  \ref{gluing-triangles} also $T=\T(xyz)$ does.
  \\
  \noindent{}
  4. Assume now that $0$ lies in the interior of $T$ (of course a
  non-degenerate triangle is a Jordan curve).  Let $\alfa: [0,L] \ra
  \Xr$ be the segment $[x,y]$ and let $ \beta_t: [0,1] \ra X$ be a
  constant speed parametrisation of the segment from $z$ to
  $\alfa(t)$.  Then
  \begin{gather*}
    F : Q=[0,1]^2 \ra X \qquad H(t,s) = \beta_t(s)
  \end{gather*}
  is a continuous map. 
Since $F\restr{\partial Q}: \partial Q \ra T $ is a degree one map
%
  $F(Q)$ must fill the interior of $T$.
  In particular there is $t_0\in (0,L)$ such that $\beta_{t_0}$ passes
  through $0$. Then we can apply the previous argument to both
  triangles $\T(xz\alfa(t_0))$ and $\T(yz\alfa(t_0))$. Applying again
  Lemma \ref{gluing-triangles} we get that the angle condition holds
  for $T$.
  \\
  \noindent{}
  5.  Finally consider the case in which $0 \not \in \overline {R}$,
  where $R$ is the interior of $T$.  Assume that
  \begin{gather}
    \label{eq:angolozzi-scaramellozzi}
    \ango \bigl (\uncobra(x), \uncobra(y) \bigr ) \geq \max \Bigl \{
    \ango \bigl (\uncobra(x), \uncobra(z) \bigr ) , \ango \bigl
    (\uncobra(z), \uncobra(y) \bigr ) \Bigr \}.
  \end{gather}
  Since $0$ does not belong to $\overline{R}$, in particular it does
  not lie on $[x,y]$. By Theorem \ref{unicita-superfinale} this
  implies
  \begin{gather*}
    \ango \bigl (\uncobra(x), \uncobra(y) \bigr ) < \pi/m.
  \end{gather*}
  Write $x=\camel(\sigma(t), s) $, $y=\camel(\sigma(t'), s') $ and
  $z=\camel(\sigma(t''), s'') $.  Then
  \eqref{eq:angolozzi-scaramellozzi} reads
  \begin{gather*}
    | \solle (t) - \solle(t')| \geq \max \bigl \{ | \solle (t) -
    \solle(t'')| , | \solle (t'') - \solle(t')|\bigr \}.
  \end{gather*}
  By interchanging $t$ and $t'$ (that is $x$ and $y$) we can then
  assume that $t\leq t'' \leq t' < \opt$.  We claim that
  $\overline{R}\subset \sect[t,t']$. Indeed if there is a point of
  $\overline{R}$ outside $\sect[t,t']$, there must be some point $w
  \in \partial R$ outside of $\sect[t,t']$.  But $\partial R = [x,y]
  \cup [y,z] \cup [x,z]$ and the three segments are contained in
  $\sect[t,t']$.
%
%
%
\end{proof}

We are now finally ready to prove the main result of the paper.

\begin{teo}  
  Let $(X,\om)$ be a \Keler curve and let $d$ be the intrinsic
  distance.
  If $\Kl$ is an upper bound for the Gaussian curvature of $ g$ on
  $\Xr$, then $(X,d)$ is a metric space of curvature $\leq \Kl$ in the
  sense of Alexandrov.
\end{teo}
\begin{proof}
  We need to prove that for any $x_0\in X$ there is a geodesic ball
  centred at $x_0$ that is a \cat--space. If $x_0\in \Xr$ this is
  well-known (Prop. \ref{AlexRiem}). If $x_0$ is an analytically
  irreducible singular point (i.e. a single branch singularity),
  thanks to Cor. \ref{localizzo} it is enough to consider the
  situation envisaged in \Sec\Sec
  \ref{sec-regu}--\ref{convex-section}. In this case the
  \cat--property of sufficiently small balls is what we have just
  proven (Theorem \ref{caso-irriducibile}). Finally we have to
  consider the case in which $x_0$ is a singular point and $X$ is
  analytically reducible at $x_0$. Let $U$ be a neighbourhood of $x$
  such that $U=U_1 \cup \cdots \cup U_N$ where $U_j$ are the
  irreducible components of $U$, $x_0\in U_j$ for each $j$ and the
  singular set of $U_j$ contains at most $x_0$.  Denote by $d_j$ the
  intrinsic distance of $(U_j, \om\restr{U_j})$.  For $r>0$ let
  $\Bd(x_0,r)$ be the geodesic ball in $(X,d)$, as usual, and let
  $\Bd_j(r)$ be the geodesic ball of radius $r$ centred at $x_0$ in
  the space $(U_j,d_j)$. By choosing $r>0$ small enough we can assume
  that any pair of points in $\Bd(x_0,r)$ is joined by a segment in
  $U$.  This follows from Lemma \ref{ball-hopf}.  It is clear that
  $\Bd_j(r) \subset \Bd(x_0,r)$.  On the other hand if $x\in
  \Bd(x_0,r)$ and $\alfa: [0,L] \ra U$ is a segment from $x_0$ to $x$ then
  $\alfa(t) \neq x_0 $ for $t>0$. So $\alfa([0,L]) \subset U_j$ for
  some $j$. Since $L(\alfa) =d(0,x) <r$ it follows that $x\in
  \Bd_j(r)$ and that $d_j(0,x) = d(0,x)$.  This shows that
  \begin{gather*}
    \Bd(x_0, r) = \Bd_1(r) \cup \cdots{} \cup \Bd_N(r).
  \end{gather*}
Moreover if $j \neq k$ any segment joining $x
  \in \Bd_j(r)$ to $ y \in \Bd_k(r)$ necessarily passes through $x_0$. 
Therefore
\begin{gather}
  \label{eq:gluglu}
  d(x,y) =  
  \begin{cases}
    d_j(x,y) & \text {if } x,   y \in \Bd_j(r) \\
    d_j(x,0) + d_k (0,y) & \text {if } x \in \Bd_j(r), y \in \Bd_k(r) ,
    j\neq k.
  \end{cases}
\end{gather}
Since each $\Bd_j(r)$ is either smooth or analytically irreducible, by
further decreasing $r$ we can assume that each $\Bd_j(r)$ is
geodesically convex in $(U_j,d_j)$ and is a \cat--space with the
distance $d_j$.  It follows from this and \eqref{eq:gluglu} that
geodesic segments are unique in $\Bd(x_0,r)$.  Let $T=\T(xyz)$ be a
triangle in $\Bd(x_0,r)$. If the three points lie in the same
$\Bd_j(r)$ the result follows from the \cat--space property of
$\Bd_j(r)$.  If $x\in  \Bd_1(r)$ and $y\in \Bd_2(r)$ and $z\in
\Bd_3(r)$, then the triangle is a tree with three edges. All the
angles vanish and the angle condition is trivially satisfied.  Finally
assume that $x,y \in \Bd_1(r)$ and $z\in \Bd_2(r)$. Then $x_0 \in
[x,z] \cap [y,z]$, so the angle at $z$ vanishes, while the angles at
$x$ and $y$ are the same as in $T'=\T(xy0)$. Since $T'\subset
\Bd_1(r)$ the angles in $T'$ are smaller than the ones in the
comparison triangle $\overline{T'}$. But $\overline {T}$ is obtained
by ``straightening'' $\overline{T'}$. Thus it follows from Alexandrov
lemma \cite[Lemma 4.3.3 p.115]{burago-burago-ivanov} that the angle
condition holds for $T$ too.
\end{proof}

The argument in the last part of the proof is the same as in
Reshetnyak Theorem \cite[p. 316]{burago-burago-ivanov}. Our case is
the simplest possible one, since the spaces are glued along a set that
consists of a single point.

\begin{teo}\label{no-low}
  If $(X,\om)$ is a \Keler curve and $x_0$ is a singular point, every
  geodesic arriving at $x_0$ branches into a continuum of different
  segments. In particular as soon as $ X_\sing \neq \vacuo$, there is
  no $\Kl \in \R$ such that $(X,d)$ be a metric space of curvature
  $\geq \Kl$ (in the sense of Alexandrov).
\end{teo}
\begin{proof}
  Assume that $ \Bd(x_0, r) = \Bd_1(r) \cup \cdots{} \cup \Bd_N(r)$ as
  above.
%
  If $N=1$ the singularity is analytically irreducible and the claim
  is already contained in Theorem \ref{unicita-superfinale}.  If $N>1$
  fix $x \in \Bd_1(r)$, $x\neq x_0$.  For any $y\in \Bd_2(r) \setminus
  \{x_0\}$ the segment from $x$ to $y$ passes through $x_0$.  This
  proves that there infinitely many segments prolonging the segment
  from $x$ to $x_0$. Since segments cannot branch in Alexandrov spaces
  with curvature bounded below it follows that no such bound can hold
  on $(X,d)$.
\end{proof}

\begin{remark}
  The point of the above result is that $\inf_{\Xr} K$ can in fact be
  finite even when $X$ contains singularities. For example consider
  $X=\{(x,y) \in \C^2: y^2=x^n\}$ with the Euclidean metric. A simple
  computation using \eqref{eq:B-uguale} shows that for $n > 4$ the
  Gaussian curvature is bounded near $(0,0)$. Nevertheless by
  Thm. \ref{no-low} there is no lower bound in the sense of
  Alexandrov.
\end{remark}
\begin{remark}
  In the case of an irreducible singularity it would be interesting to
  understand if different segments starting at the singular point can
  have the same initial tangent vector. If this were not the case the
  map $\solle$ in \eqref{eq:def-solle} would be strictly increasing
  and $\cobra\restr {C}$ would be a homeomorphism of $C$ onto
  $S^1\times \{0\}\subset C_0(X)$. Its inverse would share many
  properties of the exponential map of a Riemannian manifold. We leave
  the analysis of this problem for the future.
\end{remark}

\def\cprime{$'$}

\end{document}